\documentclass[preprint,10pt]{elsarticle}

\usepackage{lineno,hyperref}
\usepackage{amsmath}
\usepackage{subfigure}
\usepackage{caption}
\usepackage{xcolor}
\usepackage{amssymb}
\usepackage{lipsum}
\usepackage{amsfonts}
\usepackage{graphicx}
\usepackage{epstopdf}
\usepackage{algorithmic}
\usepackage{amsmath}
\usepackage{bm}
\usepackage{mathrsfs}
\usepackage{amsthm}

\newtheorem{remark}{Remark}

\modulolinenumbers[5]

\journal{Journal of Fluids and Structures}

\biboptions{numbers,sort&compress}
\begin{document}

\begin{frontmatter}

\title{A theoretical and experimental investigation of a family of immersed finite element methods} %\tnoteref{mytitlenote}


\author[]{Yongxing Wang\corref{mycorrespondingauthor}}
\cortext[mycorrespondingauthor]{Corresponding author}
\ead{scsywan@leeds.ac.uk/yongxingwang6@gmail.com}

\author[]{Peter K. Jimack}
\author[]{Mark A. Walkley}

\address{School of Computing, University of Leeds, Leeds, UK, LS2 9JT}

\begin{abstract}
In this article we consider the widely used immersed finite element method (IFEM), in both explicit and implicit form, and its relationship to our more recent one-field fictitious domain method (FDM). We review and extend the formulation of these methods, based upon an operator splitting scheme, in order to demonstrate that both the explicit IFEM and the one-field FDM can be regarded as particular linearizations  of the fully implicit IFEM. However, the one-field FDM can be shown to be more robust than the explicit IFEM and can simulate a wider range of solid parameters with a relatively large time step. In addition, it can produce results almost identical to the implicit IFEM but without iteration inside each time step. We study the effect on these methods of variations in viscosity and density of fluid and solid materials. The advantages of the one-field FDM within the IFEM framework are illustrated through a selection of parameter sets for two benchmark cases.
\end{abstract}


\end{frontmatter}

%\linenumbers

%%%%%%%%%%%%%%%%%%%%%%%%%%%%%%%%%%%%
\section{Introduction}
\label{sec_introduction}
Three major questions arise when considering a finite element method for the problem of Fluid-Structure Interactions (FSI): (1) what kind of meshes are used (interface fitted or unfitted); (2) how to couple the fluid-structure interactions (monolithic/fully-coupled or partitioned/segregated); (3) what variables are solved (velocity and/or displacement). Combinations of the answers to these questions lead to different types of numerical method. For example, \cite{Degroote_2009, K_ttler_2008}  solve for fluid velocity and solid displacement sequentially (partitioned/segregated) using an Arbitrary Lagrangian-Eulerian (ALE) fitted mesh, whereas \cite{Heil_2004, Heil_2008, Muddle_2012} use an ALE fitted mesh to solve for fluid velocity and solid displacement simultaneously (monolithic/fully-coupled) with a Lagrange Multiplier to enforce the continuity of velocity/displacement on the interface. The Immersed Finite Element Method (IFEM) \cite{peskin2002immersed, Wang_2011, Wang_2009,Wang_2013,Zhang_2007,zhang2004immersed} and the Fictitious Domain Method (FDM) \cite{baaijens2001fictitious,Boffi_2016,Glowinski_2001,Hesch_2014,Kadapa_2016,Yu_2005} use two meshes to represent the fluid and solid separately. Although IFEM could be monolithic \cite{Boffi_2015}, the classical IFEM only solves for velocity, while the solid information is arranged on the right-hand side of the fluid equation as a prescribed force term. Conversely, although the FDM may be partitioned \cite{Yu_2005}, usually the FDM approach solves for both velocity in the whole domain (fluid plus solid) and displacement of the solid simultaneously via a distributed Lagrange multiplier (DLM) to enforce the consistency of velocity/displacement in the overlapped solid domain. In the case of one-field and monolithic numerical methods for FSI problems, \cite{Auricchio_2014} introduces a 1D model using a one-field FDM formulation based on two meshes, and \cite{Hecht_2017,Pironneau_2016} introduces an energy stable monolithic method (in 2D) based on one Eulerian mesh and discrete remeshing. 

In a previous study \cite{Wang_2017}, we present a one-field monolithic fictitious domain method (subsequently referred to as the one-field FDM) which has the following main features: (1) only one velocity field is solved in the whole domain, based upon the use of an appropriate $L^2$ projection; (2) the fluid and solid equations are solved monolithically. The primary purpose of this paper is to highlight the relationship between the one-field FDM and the IFEM family of methods: demonstrating that the former is in fact a new variant of the latter which possesses a number of practical advantages. Before describing these in detail however we briefly provide further context for the one-field FDM based on comparing its features with those of existing monolithic schemes. 

FDM/DLM methods \cite{baaijens2001fictitious,Boffi_2016,Glowinski_2001,Hesch_2014,Kadapa_2016,Yu_2005} solve the solid equation, but for a displacement field, and couple this displacement with the velocity of the fictitious fluid via a Lagrange multiplier. This leads to a large discrete linear algebra system. The one-field FDM solves a smaller discrete system since it rewrites the solid equation in terms of a velocity variable and couples the fictitious fluid through a finite element interpolation. Monolithic Eulerian methods \cite{Hecht_2017, Pironneau_2016} also express the solid equation in terms of velocity, and the fluid and solid are coupled naturally on an interface-fitted mesh. The one-field FDM also uses two meshes to represent the fluid and solid respectively. Consequently, before discretization in space, these two methods have many similarities, the advantage of the one-field FDM being that interface fitting is not required.

As discussed above, the primary purpose of this paper is to demonstrate that the proposed one-field FDM is a particular linearization of the fully implicit IFEM. It is more robust than the classical explicit IFEM, and presents almost the same accuracy as the implicit IFEM. We will show: 
\begin{itemize}
\item [(1)] The one-field FDM is based upon the implicit expression of the FSI force. This FSI force term is linearized using the velocity at the current configuration (instead of displacement at the reference configuration), which is an approximation of the fully implicit IFEM \cite{wang2006immersed, wang2007iterative, wang2009computational} but without requiring a nonlinear iteration at each time step.
\item [(2)] In the simple case of equal density and viscosity for both fluid and solid, the only difference between the one-field FDM and the explicit IFEM (explicitly expressing the FSI force term) \cite{Wang_2011, Wang_2009,Wang_2013,Zhang_2007,zhang2004immersed} is that there are some terms of order $O(\Delta t)$ and $O(\Delta t^2)$ retained in the former. However these terms have a helpful stabilizing effect, which can allow the one-field FDM to use a time step that is significantly larger than the explicit IFEM, with almost the same accuracy.
\item [(3)] The one-field FDM can naturally deal with the case of different densities and different viscosities in the fluid and solid.
\end{itemize}

The paper is organized as follows. The control equations and a general finite element weak formulation are introduced in Section \ref{control_equations} and \ref{weak_formulation} respectively, followed by a dimensionless weak formulation in Section \ref{dimensionless_weak} and time discretization in Section \ref{time_discretization}. Different types of IFEM approaches are then discussed in Section \ref{diferent_ifem}, and their relationship with the weak form of Section \ref{dimensionless_weak} is highlighted. An operator splitting scheme is introduced in Section \ref{splitting_scheme} followed by a comparison between the IFEM and one-field FDM approaches in Section \ref{ifem} and \ref{one_field_fdm}. Numerical examples are given in Section \ref{sec_numerical_exs}, and conclusions are presented in Section \ref{sec_conclusions}.

%%%%%%%%%%%%%%%%%%%%%%%%%%%%%%%%%%%%%%
\section{Control equations}
\label{control_equations}
\begin{figure}[h!]
\centering
\includegraphics[width=2.5in,angle=0]{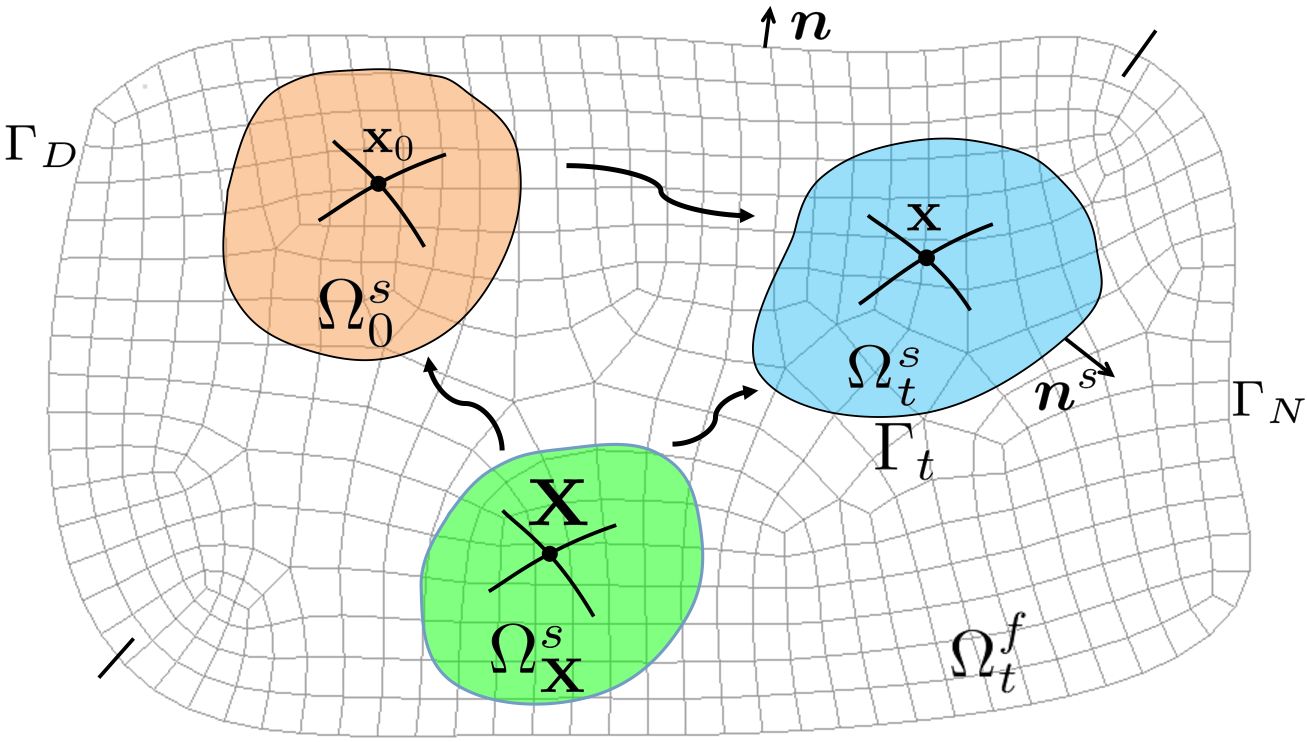}
\caption {\scriptsize Schematic diagram of FSI, $\Omega=\Omega_t^f\cup \Omega_t^s$ and $\Gamma=\Gamma_D\cup \Gamma_N$.} 
\label{fig1:Schematic diagram of FSI}
\end{figure}

In the following context, $\Omega_t^f\subset\mathbb{R}^d$ and $\Omega_t^s\subset\mathbb{R}^d$ with $d=2,3$ denote the fluid and solid domain respectively which are time dependent regions as shown in Figure \ref{fig1:Schematic diagram of FSI}. $\Omega=\Omega_t^f \cup \Omega_t^s $ is a fixed domain (with outer boundary $\Gamma$) and $\Gamma_t=\partial\Omega_t^f\cap\partial\Omega_t^s$ is the moving interface between fluid and solid. We denote by ${\bf X}$ the reference (material) coordinates of the solid, by ${\bf x}={\bf x}(\cdot,t)$ the current coordinates of the solid, and by ${\bf x}_0$ the initial coordinates of the solid.

Let $\rho, \mu, {\bf u}, {\bm{\sigma}}$ and ${\bf g}$ denote the density, viscosity, velocity, stress tensor and acceleration due to gravity respectively. We assume both an incompressible fluid in $\Omega_t^f$ and incompressible solid in $\Omega_t^s$. The conservation of momentum and conservation of mass take the same form in the fluid and solid (just differing in the specific expressions of ${\bm\sigma}$):

Momentum equation:
\begin{equation} \label{momentum_equation}
\rho\frac{d{\bf u}}{dt}
=\nabla \cdot {\bm\sigma}+\rho{\bf g},
\end{equation}
and continuity equation:
\begin{equation} \label{continuity_equation}
\nabla \cdot {\bf u}=0.
\end{equation}

An incompressible Newtonian constitutive equation in $\Omega_t^f$ can be expressed as:
\begin{equation} \label{constitutive_fluid}
{\bm\sigma}={\bm\sigma}^f={\bm\tau}^f-p^f{\bf I},
\end{equation}
where  ${\rm D}{\bf u}=\nabla {\bf u}+\nabla^{\scriptsize T} {\bf u}$, and
\begin{equation} \label{constitutive_fluid_deviatoric}
{\bm\tau}^f=\mu^f {\rm D}{\bf u}^f
\end{equation}
is the deviatoric part of stress ${\bm\sigma}^f$. 

In $\Omega_t^s$ we use an incompressible neo-Hookean solid  with viscosity $\mu^s$ \cite{Boffi_2016,Yu_2005,baaijens2001fictitious} (see appendices for a brief discussion of other solid material models). The constitutive equation may be expressed as:
\begin{equation} \label{constitutive_solid}
{\bm\sigma}={\bm\sigma}^s={\bm\tau}^s-p^s{\bf I},
\end{equation}
with
\begin{equation} \label{constitutive_solid_deviatoric}
{\bm\tau}^s=c_1\left({\bf F}{\bf F}^T-{\bf I}\right)+\mu^s{\rm D}{\bf u}^s
\end{equation}
being the deviatoric part of stress ${\bm\sigma}^s$, where 
\begin{equation}\label{definition_f}
{\bf F}
=\frac{\partial {\bf x}}{\partial {\bf X}}
=\frac{\partial {\bf x}}{\partial {\bf x}_0} \frac{\partial {\bf x}_0}{\partial {\bf X}}
=\nabla_0{\bf x}\nabla_{\bf X}{\bf x}_0
\end{equation}
is the deformation tensor of the solid, and $c_1$ a solid material parameter. 

Finally the system is completed with the following boundary and initial conditions.

Interface continuity conditions:
\begin{equation}\label{interfaceBC1}
{\bf u}^f={\bf u}^s\quad on \quad  \Gamma_t
\end{equation}
and
\begin{equation}\label{interfaceBC2}
{\bm \sigma}^f{\bf n}^s= {\bm \sigma}^s{\bf n}^s\quad on \quad  \Gamma_t.
\end{equation}
Dirichlet and Neumann boundary conditions:
\begin{equation}\label{homogeneous_boundary}
{\bf u}^f=\bar{\bf u}\quad on \quad  \Gamma_D
\end{equation}
and
\begin{equation}\label{Neumann}
{\bm \sigma}^f{\bf n}^s= \bar{\bf h}\quad on \quad  \Gamma_N,
\end{equation}
where $\Gamma=\Gamma_D\cup\Gamma_N$. Initial conditions:
\begin{equation} \label{initialcd_fluid}
\left. {\bf u}^f\right|_{t=0}={\bf u}_0^f
\end{equation}
and
\begin{equation} \label{initialcd_solid}
\left. {\bf u}^s\right|_{t=0}={\bf u}_0^s.
\end{equation}

%%%%%%%%%%%%%%%%%%%%%%%%%%%%%%%%%%%%%%%%%%%%%%%%%%%%%%%%%
\section{Weak formulation}
\label{weak_formulation}
In the following context, let $L^2(\omega)$ be the square integrable functions in domain $\omega$, endowed with norm $\left\|u\right\|_{0,\omega}^2=\int_\omega \left|u\right|^2$ ($u\in L^2(\omega)$). Let $H^1(\omega)=\left\{u:u, \nabla u\in L^2(\omega)\right\}$ with the norm denoted by $\left\|u\right\|_{1,\omega}^2=\left\|u\right\|_{0,\omega}^2+\left\|\nabla u\right\|_{0,\omega}^2$. We also denote by $H_0^1(\omega)$ the subspace of $H^1(\omega)$ whose functions have zero values on the Dirichlet boundary of $\omega$.

Let
$
{p}=\left \{ 
\begin{matrix}
{{p}^f \quad in \quad \Omega_t^f} \\
{{p}^s \quad in \quad \Omega_t^s} \\
\end{matrix}\right.
$.
Given ${\bf v}\in H_0^1(\Omega)^d$, we perform the following symbolic operations:
\begin{equation*}
\begin{split}
\int_{\Omega}{\rm Eq.}(\ref{momentum_equation})\left({\bm\sigma}\right)\cdot{\bf v}
&\equiv
\int_{\Omega_t^f}{\rm Eq.}(\ref{momentum_equation})\left({\bm\sigma}^f\right)\cdot{\bf v}
+\int_{\Omega_t^s}{\rm Eq.}(\ref{momentum_equation})\left({\bm\sigma}^s\right)\cdot{\bf v}\\
&\equiv
\int_{\Omega}{\rm Eq.}(\ref{momentum_equation})\left({\bm\sigma}^f\right)\cdot{\bf v}
+\int_{\Omega_t^s}\left({\rm Eq.}(\ref{momentum_equation})\left({\bm\sigma}^s\right)-{\rm Eq.}(\ref{momentum_equation})\left({\bm\sigma}^f\right)\right)\cdot{\bf v}.
\end{split}
\end{equation*}

Integrating the stress terms by parts, the above operations, using constitutive equation (\ref{constitutive_fluid}) and (\ref{constitutive_solid}) and interface condition (\ref{interfaceBC2}), give:
\begin{equation}\label{operation1}
\begin{split}
&\rho^f\int_{\Omega}\frac{d{\bf u}}{dt} \cdot{\bf v}
+\int_{\Omega}{\bm{\tau}^f}:\nabla{\bf v}
-\int_{\Omega}p\nabla \cdot {\bf v} \\
+&\left(\rho^s-\rho^f\right)\int_{\Omega_t^s}\frac{d{\bf u}}{dt}\cdot{\bf v}
+\int_{\Omega_t^s}\left({\bm{\tau}^s}-{\bm{\tau}^f}\right):\nabla{\bf v} \\
=&\int_{\Omega}\rho^f{\bf g}\cdot{\bf v}
+\int_{\Omega_t^s}\left(\rho^s-\rho^f\right){\bf g}\cdot{\bf v}
+\int_{\Gamma_N}\bar{\bf h}\cdot{\bf v}.
\end{split}
\end{equation}
Note that the integrals on the interface $\Gamma_t$ are cancelled out using boundary condition (\ref{interfaceBC2}), because they are internal forces for the whole FSI system. Combining with the following symbolic operations for $q\in L^2(\Omega)$,
$$
-\int_{\Omega_t^f}{\rm Eq.}(\ref{continuity_equation})q
-\int_{\Omega_t^s}{\rm Eq.}(\ref{continuity_equation})q
\equiv
-\int_{\Omega}{\rm Eq.}(\ref{continuity_equation})q,
$$ 
leads to the weak form of the FSI system as follows.

Given ${\bf u}_0$ and $\Omega_0^s$, for each $t>0$ find ${\bf u}(t)\in H^1(\Omega)^d$, $p(t) \in L^2(\Omega)$ and $\Omega_t^s$, such that $\forall {\bf v}\in H_0^1(\Omega)^d$, $\forall q \in L^2(\Omega)$, the following two equations hold:
\begin{equation}\label{weak_form1}
\begin{split}
&\rho^f\int_{\Omega}\frac{\partial{\bf u}}{\partial t} \cdot{\bf v}
+\rho^f \int_{\Omega}\left({\bf u}\cdot\nabla\right){\bf u}\cdot{\bf v}
+\frac{\mu^f}{2}\int_{\Omega}{\rm D}{\bf u}:{\rm D}{\bf v}
-\int_{\Omega}p\nabla \cdot {\bf v} \\
+&\rho^{\delta}\int_{\Omega_t^s}  \frac{\mathfrak{d}{\bf u}}{\mathfrak{d} t}    \cdot{\bf v}
+\frac{\mu^\delta}{2}\int_{\Omega_t^s}{\rm D}{\bf u}:{\rm D}{\bf v}
+c_1\int_{\Omega_t^s}\left({{\bf F}{\bf F}^T-{\bf I}}\right):\nabla{\bf v} \\
=&\int_{\Omega}\rho^f{\bf g}\cdot{\bf v}
+\int_{\Omega_t^s}\rho^\delta{\bf g}\cdot{\bf v}
+\int_{\Gamma_N}\bar{\bf h}\cdot{\bf v},
\end{split}
\end{equation}
 (where $\rho^\delta=\rho^s-\rho^f$ and $\mu^\delta=\mu^s-\mu^f$) and
\begin{equation}\label{weak_form2}
-\int_{\Omega} q\nabla \cdot {\bf u}=0.
\end{equation}		

In equation (\ref{weak_form1}), $\frac{\mathfrak{d}}{\mathfrak{d}t}$ is the time derivative with respect to a frame moving with the solid velocity ${\bf u}^s=\left. {\bf u}\right|_{\Omega_t^s}$.

%%%%%%%%%%%%%%%%%%%%%%%%%%%%%%%%%%%%%%%%%%%%%%%%%%%%%%%%%
\section{Dimensionless weak formulation}
\label{dimensionless_weak}
We may rewrite the weak form (\ref{weak_form1}) and (\ref{weak_form2}) in a dimensionless form by introducing the following scales \cite{Yu_2005}: $\tilde{L}$ for length, $\tilde{U}$ for velocity, $\tilde{L}/\tilde{U}$ for time and $\rho^f\tilde{U}^2$ for pressure $p$ and parameter $c_1$. Dividing by $\rho^f\tilde{U}^3/\tilde{L}$ on both sides of equation (\ref{weak_form1}), we have:
\begin{equation}\label{weak_form1_dimensionless}
\begin{split}
&\int_{\Omega}\frac{\partial\tilde{\bf u}}{\partial \tilde {t}} \cdot\tilde{\bf v}
+\int_{\Omega}\left(\tilde{\bf u}\cdot\nabla\right)\tilde{\bf u}\cdot\tilde{\bf v}
+\frac{1}{2Re}\int_{\Omega}{\rm D}\tilde{\bf u}:{\rm D}\tilde{\bf v}
-\int_{\Omega}\tilde{p}\nabla \cdot \tilde{\bf v} \\
+&\left(\rho^{r}-1\right)\int_{\Omega_t^s} \frac{\mathfrak{d}\tilde{\bf u}}{\mathfrak{d} \tilde{t}}
 \cdot\tilde{\bf v}
+\frac{\mu^r-1}{2Re}\int_{\Omega_t^s}{\rm D}\tilde{\bf u}:{\rm D}\tilde{\bf v}
+\tilde{c}_1\int_{\Omega_t^s}\left({{\bf F}{\bf F}^T-{\bf I}}\right):\nabla{\bf v} \\
=& Fr\int_{\Omega}\frac{{\bf g}}{|{\bf g}|}\cdot\tilde{\bf v}
+Fr\left(\rho^r-1\right)\int_{\Omega_t^s}\frac{{\bf g}}{|{\bf g}|}\cdot\tilde{\bf v}
+\int_{\Gamma_N}\tilde{\bf h}\cdot\tilde{\bf v},
\end{split}
\end{equation}
where $\tilde{\bf u}={\bf u}/\tilde{U}$, $\tilde{\bf v}={\bf v}/\tilde{U}$, $\tilde{p}=p/\rho^f\tilde{U}^2$, $\tilde{t}=t\tilde{U}/\tilde{L}$, $\tilde{\bf h}=\bar{\bf h}\tilde{L}/\rho^f\tilde{U}^2$ and the following parameters.
\begin{itemize}
\item Reynolds number:
\begin{equation}
Re=\rho^f\tilde{U}\tilde{L}/\mu^f.
\end{equation}
\item density ratio:
\begin{equation}
\rho^r=\rho^s/\rho^f.
\end{equation}
\item viscosity ratio:
\begin{equation}
\mu^r=\mu^s/\mu^f.
\end{equation}
\item material parameter:
\begin{equation}
\tilde{c}_1=c_1/\rho^f\tilde{U}^2.
\end{equation}
\item Froude number:
\begin{equation}\label{froude}
Fr=|{\bf g}|\tilde{L}/\tilde{U}^2.
\end{equation}
\end{itemize}
The dimensionless form of equation (\ref{weak_form2}) is straightforward by dividing by $\rho^f\tilde{U}^3/\tilde{L}$:
\begin{equation}\label{weak_form2_dimensionless}
-\int_{\Omega} \tilde{q}\nabla \cdot \tilde{\bf u}=0,
\end{equation}
with $\tilde{q}=q/\rho^f\tilde{U}^2$. For the sake of notation convenience, we shall still use ${\bf u}$, ${\bf v}$, $p$ and $q$ instead of $\tilde{\bf u}$, $\tilde{\bf v}$, $\tilde{p}$ and $\tilde{q}$ in equation (\ref{weak_form1_dimensionless}) and (\ref{weak_form2_dimensionless}) in the following context.

%%%%%%%%%%%%%%%%%%%%%%%%%%%%%%%%%%%%%%%%%%%%%%%%%%%%%%%%
\section{Discretization in time}
\label{time_discretization}
Using the backward Euler method to discretize in time, equations (\ref{weak_form1_dimensionless}) and (\ref{weak_form2_dimensionless}) may be discretized as follows:

Given ${\bf u}_n$, $p_n$ and $\Omega_n^s$, find ${\bf u}_{n+1}\in H^1(\Omega)^d$, $p_{n+1} \in L^2(\Omega)$ and $\Omega_{n+1}^s$, such that $\forall{\bf v}\in H_0^1\left(\Omega\right)^d$, $\forall q\in L^2(\Omega)$, the following two equations hold:
\begin{equation}\label{weak_form1_time}
\begin{split}
&\int_{\Omega}\frac{{\bf u}_{n+1}-{\bf u}_n}{\Delta t} \cdot{\bf v}
+\int_{\Omega}\left({\bf u}_{n+1}\cdot\nabla\right){\bf u}_{n+1}\cdot{\bf v} \\
+&\frac{1}{2Re}\int_{\Omega}{\rm D}{\bf u}_{n+1}:{\rm D}{\bf v}
-\int_{\Omega}p_{n+1}\nabla \cdot {\bf v}\\
+&\left(\rho^{r}-1\right)\int_{\Omega_{n+1}^s}\frac{{\bf u}_{n+1}-{\bf u}_n}{\Delta t} \cdot{\bf v}
+\frac{\mu^r-1}{2Re}\int_{\Omega_{n+1}^s}{\rm D}{\bf u}_{n+1}:{\rm D}{\bf v} \\
+& \tilde{c}_1\int_{\Omega_{n+1}^s}\left({{\bf F}{\bf F}^T-{\bf I}}\right):\nabla{\bf v} \\
=& Fr\int_{\Omega}\frac{{\bf g}}{|{\bf g}|}\cdot{\bf v}
+Fr\left(\rho^r-1\right)\int_{\Omega_{n+1}^s}\frac{{\bf g}}{|{\bf g}|}\cdot{\bf v} 
+\int_{\Gamma_N}\tilde{\bf h}\cdot{\bf v},
\end{split}
\end{equation}
and
\begin{equation}\label{weak_form2_time}
-\int_{\Omega} q\nabla \cdot {\bf u}_{n+1}=0.
\end{equation}

\begin{remark}
$\Omega_{n+1}^s$ is updated from $\Omega_n^s$ by the following formula:
\begin{equation}\label{update_of_solid_mesh}
\Omega_{n+1}^s=\left\{{\bf x}:{\bf x}={\bf x}_n+\Delta t{\bf u}_{n+1}, {\bf x}_n\in\Omega_n^s \right\}.
\end{equation}
\end{remark}

\section{Implementations of different IFEM approaches}
\label{diferent_ifem}
Having introduced the weak formulation and time discretization in the previous sections, we now consider implementation details, and demonstrate that different choices lead to methods that correspond to existing IFEM schemes, as well as the proposed one-field FDM \cite{Wang_2017}. We can see from (\ref{weak_form1_time}) that the integrals are carried out over two different domains: the whole domain $\Omega$ and the moving solid domain $\Omega_{n+1}^s$. The IFEM methods compute these two types of integrals based on two different meshes, and use projection/distribution functions to transfer information between the two meshes \cite{baaijens2001fictitious,Boffi_2016,Glowinski_2001,Hesch_2014,Kadapa_2016,Yu_2005,Wang_2011, Wang_2009,Wang_2013,Zhang_2007,zhang2004immersed}. The one-field FDM also adopts two meshes, and uses the FEM interpolation function (also used in \cite{Wang_2009}) to transfer information between the two meshes. In the remainder of this section we focus on how different IFEM approaches fit within the framework of this weak formulation.

The classical IFEM is introduced in \cite{zhang2004immersed, Zhang_2007}, in which all the solid integrals (in $\Omega_{n+1}^s$) are moved to the right-hand side of equation (\ref{weak_form1_time}) as a force term and evaluated at the previous time step as follows:
\begin{equation}\label{weak_form_classical_ifem}
\begin{split}
&\int_{\Omega}\frac{{\bf u}_{n+1}-{\bf u}_n}{\Delta t} \cdot{\bf v}
+\int_{\Omega}\left({\bf u}_{n+1}\cdot\nabla\right){\bf u}_{n+1}\cdot{\bf v} \\
+&\frac{1}{2Re}\int_{\Omega}{\rm D}{\bf u}_{n+1}:{\rm D}{\bf v}
-\int_{\Omega}p_{n+1}\nabla \cdot {\bf v}\\
=&\left(1-\rho^{r}\right)\int_{\Omega_n^s}\frac{{\bf u}_{n}-{\bf u}_{n-1}}{\Delta t} \cdot{\bf v}
-\boxed{\frac{\mu^r-1}{2Re}\int_{\Omega_n^s}{\rm D}{\bf u}_n:{\rm D}{\bf v}} \\
-& \tilde{c}_1\int_{\Omega_n^s}\left({{\bf F}_n{\bf F}_n^T-{\bf I}}\right):\nabla{\bf v} \\
+& Fr\int_{\Omega}\frac{{\bf g}}{|{\bf g}|}\cdot{\bf v}
+Fr\left(\rho^r-1\right)\int_{\Omega_n^s}\frac{{\bf g}}{|{\bf g}|}\cdot{\bf v} 
+\int_{\Gamma_N}\tilde{\bf h}\cdot{\bf v}.
\end{split}
\end{equation}
The above formulation differs from \cite{zhang2004immersed, Zhang_2007} only in the following respects: 
\begin{enumerate}
\item [(1)]The boxed term in (\ref{weak_form_classical_ifem}) vanishes in \cite{zhang2004immersed, Zhang_2007} because the fluid stress within the solid domain is neglected (which is equivalent to setting $\mu^r=1$).
\item [(2)] \cite{zhang2004immersed, Zhang_2007} use the stabilized equal-order finite element method to treat convection and pressure after discretization in space, while we shall use a splitting scheme to deal with convection and a stable mixed-order finite element space for the velocity and pressure.
\item [(3)] The above formulation is expressed in a dimensionless form.
\end{enumerate}
However these differences are not the distinguishing features of IFEM, and do not influence any comparisons (in Section \ref{sec_numerical_exs} we show that our implementation of IFEM produces the same results as in the literature). The distinguishing feature of IFEM is its development from the Immersed Boundary Method (IBM) \cite{peskin2002immersed}: the solid information is based on the previous time step and arranged on the right-hand side of the fluid equation as a force term, which is computed on the solid mesh, distributed to the fluid mesh and then added to the fluid equation as an extra term.

The IFEM formulation (\ref{weak_form_classical_ifem}) represents an explicit forcing strategy, which approximates the time derivative in the solid using values from the previous two time steps. Errors may be accumulated as time evolves in this case, and this may be observed in numerical tests (see, for example, Figures \ref{cavity_solid_deformation_parameter3} and \ref{cavity_solid_deformation_parameter4}). There is a semi-implicit formulation which introduces an indicator function $I({\bf x})$ to smear out the densities across the fluid-solid interface \cite{Wang_2011}. Based on this indicator function $I({\bf x})$ (see \cite{Wang_2011} for the definition), the formulation (\ref{weak_form_classical_ifem}) may be modified as follows:
\begin{equation}\label{weak_form_classical_semi_ifem}
\begin{split}
&\left(1+\left(\rho^{r}-1\right)I\left(\bf x\right)\right)\int_{\Omega}\frac{{\bf u}_{n+1}-{\bf u}_n}{\Delta t} \cdot{\bf v}
+\int_{\Omega}\left({\bf u}_{n+1}\cdot\nabla\right){\bf u}_{n+1}\cdot{\bf v} \\
+&\frac{1}{2Re}\int_{\Omega}{\rm D}{\bf u}_{n+1}:{\rm D}{\bf v}
-\int_{\Omega}p_{n+1}\nabla \cdot {\bf v}\\
=&\frac{1-\mu^r}{2Re}\int_{\Omega_n^s}{\rm D}{\bf u}_n:{\rm D}{\bf v}
-\tilde{c}_1\int_{\Omega_n^s}\left({{\bf F}_n{\bf F}_n^T-{\bf I}}\right):\nabla{\bf v} \\
+& Fr\int_{\Omega}\frac{{\bf g}}{|{\bf g}|}\cdot{\bf v}
+Fr\left(\rho^r-1\right)\int_{\Omega_n^s}\frac{{\bf g}}{|{\bf g}|}\cdot{\bf v} 
+\int_{\Gamma_N}\tilde{\bf h}\cdot{\bf v}.
\end{split}
\end{equation}
Furthermore, a fully implicit forcing strategy may also be considered as follows:
\begin{equation}\label{weak_form_classical_implicit_ifem}
\begin{split}
&\int_{\Omega}\frac{{\bf u}_{n+1}-{\bf u}_n}{\Delta t} \cdot{\bf v}
+\int_{\Omega}\left({\bf u}_{n+1}\cdot\nabla\right){\bf u}_{n+1}\cdot{\bf v} \\
+&\frac{1}{2Re}\int_{\Omega}{\rm D}{\bf u}_{n+1}:{\rm D}{\bf v}
-\int_{\Omega}p_{n+1}\nabla \cdot {\bf v}\\
=&\left(1-\rho^{r}\right)\int_{\Omega_n^s}\frac{{\bf u}_{n+1}-{\bf u}_n}{\Delta t} \cdot{\bf v}
-\frac{\mu^r-1}{2Re}\int_{\Omega_n^s}{\rm D}{\bf u}_{n+1}:{\rm D}{\bf v} \\
-& \tilde{c}_1\int_{\Omega_n^s}\left({{\bf F}_{n+1}{\bf F}_{n+1}^T-{\bf I}}\right):\nabla{\bf v} \\
+& Fr\int_{\Omega}\frac{{\bf g}}{|{\bf g}|}\cdot{\bf v}
+Fr\left(\rho^r-1\right)\int_{\Omega_n^s}\frac{{\bf g}}{|{\bf g}|}\cdot{\bf v} 
+\int_{\Gamma_N}\tilde{\bf h}\cdot{\bf v}.
\end{split}
\end{equation}
It can be seen that the force term on the right-hand side of above equation is computed using the velocity at the current time step, which then needs to be iteratively constructed. This idea of implicit forcing was first utilized in the immersed boundary method (IBM) \cite{newren2007unconditionally,newren2008comparison}, and then also adopted within IFEM in \cite{wang2006immersed, wang2007iterative, wang2009computational}, where a Newton-Krylov iterative procedure is used to solve this nonlinear system. In our implementation, for simplicity, we use fixed point iteration to demonstrate the implicit IFEM. The fixed point iteration generally converges more slowly than the Newton-Krylov method, however it is not our purpose to compare the efficiency of these two implicit forcing strategies. Instead we shall demonstrate that both the implicit IFEM (iterating at each time step) and the one-field FDM (which needs no iteration inside the time step) produce almost identical results.

Based upon the above discussion, we next introduce an operator splitting scheme in Section \ref{splitting_scheme}. In Section \ref{ifem} we then present an explicit and an implicit forcing strategy for IFEM (corresponding to (\ref{weak_form_classical_ifem}) and (\ref{weak_form_classical_implicit_ifem}) respectively), and in Section \ref{one_field_fdm} the one-field FDM formulation is illustrated in detail with the proposed operator splitting scheme.
%%%%%%%%%%%%%%%%%%%%%%%%%%%
\section{An operator splitting scheme}
\label{splitting_scheme}
The operator splitting scheme is introduced here in order to treat the non-linear convection term in the Navier-Stokes equation, and simplify the saddle-point problem, which is widely adopted to solve pure fluid equations \cite{Glowinski_2003, laval1990fractional}. The fluid-structure coupling process can still be designed either explicitly or implicitly inside the diffusion step (as discussed in Section \ref{ifem} and \ref{one_field_fdm}). In this article, we focus on studying the FSI coupling strategies rather than different methods to deal with the convection or saddle-point problem.
\begin{itemize}
\item[(1)] Convection step:
\begin{equation}\label{convection_step}
\int_{\Omega}\frac{{\bf u}_{n+1/3}-{\bf u}_n}{\Delta t} \cdot{\bf v}
+\int_{\Omega}\left({\bf u}_{n+1/3}\cdot\nabla\right){\bf u}_{n+1/3}\cdot{\bf v}
=0.
\end{equation}
\item[(2)] Diffusion step:
\begin{equation}\label{diffusion}
\begin{split}
&\int_{\Omega}\frac{{\bf u}_{n+2/3}-{\bf u}_{n+1/3}}{\Delta t} \cdot{\bf v}
+\frac{1}{2Re}\int_{\Omega}{\rm D}{\bf u}_{n+2/3}:{\rm D}{\bf v}  \\
+&\left(\rho^{r}-1\right)\int_{\Omega_n^s}\frac{{\bf u}_{n+2/3}-{\bf u}_n}{\Delta t} \cdot{\bf v}
+\frac{\mu^r-1}{2Re}\int_{\Omega_n^s}{\rm D}_n{\bf u}_{n+2/3}:{\rm D}_n{\bf v} \\
+& \tilde{c}_1\int_{\Omega_n^s}\left({{\bf F}_{n+2/3}{\bf F}_{n+2/3}^T-{\bf I}}\right):\nabla_n{\bf v} \\
=& Fr\int_{\Omega}\frac{{\bf g}}{|{\bf g}|}\cdot{\bf v}
+Fr\left(\rho^r-1\right)\int_{\Omega_n^s}\frac{{\bf g}}{|{\bf g}|}\cdot{\bf v} 
+\int_{\Gamma_N}\tilde{\bf h}\cdot{\bf v}.
\end{split}
\end{equation}
In the above, $\nabla_n(\cdot)=\frac{\partial(\cdot)}{\partial{\bf x}_n}$ and ${\rm D}_n=\nabla_n+\nabla_n^T$.
\item[(3)] Pressure step:
\begin{equation}\label{pressure1}
\int_{\Omega}\frac{{\bf u}_{n+1}-{\bf u}_{n+2/3}}{\Delta t} \cdot{\bf v}
-\int_{\Omega}p_{n+1}\nabla \cdot {\bf v}
=0.
\end{equation}
and
\begin{equation}\label{pressure2}
-\int_{\Omega} q\nabla \cdot {\bf u}_{n+1}=0.
\end{equation}
\end{itemize}

\begin{remark}\label{fraction_f}
Notice that the variables ${\bf u}_{n+1/3}$ or ${\bf u}_{n+2/3}$ are just intermediate values, not the velocity at time $t=t_n+\frac{\Delta t}{3}$ or $t=t_n+\frac{2\Delta t}{3}$. The notation ${\bf F}_{n+1/3}$ or ${\bf F}_{n+2/3}$ is interpreted as follows: 
\begin{equation}\label{fraction_f_eq}
{\bf F}_t=\frac{\partial {\bf x}_t}{\partial {\bf X}}
=\nabla_{\bf X}\left({\bf x}_n+{\bf u}_t\Delta t\right),
\end{equation}
with $t=n+1/3$ or $n+2/3$.
\end{remark}

\begin{remark}
Compared with the 2-step splitting scheme in our original paper \cite{Wang_2017}, this 3-step splitting scheme decouples the FSI problem and Stokes/saddle-point problem into two separate steps. The fluid and solid are coupled in the diffusion step, which may be effectively solved by the preconditioned Conjugate Gradient algorithm. The pressure step becomes a ``degenerate'' Stokes problem \cite{Glowinski_2003}, which can also be efficiently solved (readers may refer to \cite [Section 34] {Glowinski_2003} for more discussion). There are a variety of numerical methods to treat the convection equation (\ref{convection_step}), such as wave-like methods \cite{Glowinski_2003}, characteristic based schemes \cite{Glowinski_2003, Zienkiewic2014, Hecht_2017}, upwind schemes (including the Streamline Upwind Petrov Galerkin (SUPG) method) \cite{Glowinski_2003, Zienkiewic2014} or the Least-squares method \cite{Zienkiewic2014}. In our implementations we primarily use this latter approach.
\end{remark}

It can be seen that the fluid-structure interaction only occurs in the diffusion step (\ref{diffusion}) based upon this operator splitting scheme. In order to solve equation (\ref{diffusion}) the one-field FDM and IFEM use different strategies. In the following two sections we focus on this diffusion step, and present the differences between the one-field FDM, and IFEM strategies.

%%%%%%%%%%%%%%%%%%%%%%%%%%%
\section{Explicit and implicit forcing}
\label{ifem}
For notational convenience let us define the following force term:
\begin{equation}\label{force_fsi}
\begin{split}
\mathcal{F}_t^{FSI}=&\left(\rho^{r}-1\right)\int_{\Omega_n^s}\frac{{\bf u}_t-{\bf u}_n}{\Delta t} \cdot{\bf v}
+\frac{\mu^r-1}{2Re}\int_{\Omega_n^s}{\rm D}_n{\bf u}_t:{\rm D}_n{\bf v} \\
+& \tilde{c}_1\int_{\Omega_n^s}\left({{\bf F}_t{\bf F}_t^T-{\bf I}}\right):\nabla_n{\bf v},
\end{split}
\end{equation}
where $t=n+1/3$ or $n+2/3$ as in Remark \ref{fraction_f}.

Based upon the splitting scheme adopted here, we use ${\bf u}_{n+1/3}$, obtained from the previous convection step, to evaluate $\mathcal{F}_t^{FSI}$, and solve equation (\ref{diffusion}) as follows.
\begin{itemize}
\item Explicit forcing:
\begin{equation}\label{diffusion_explicit}
\begin{split}
&\int_{\Omega}\frac{{\bf u}_{n+2/3}-{\bf u}_{n+1/3}}{\Delta t} \cdot{\bf v}
+\frac{1}{2Re}\int_{\Omega}{\rm D}_n{\bf u}_{n+2/3}:{\rm D}_n{\bf v}  \\
=& Fr\int_{\Omega}\frac{{\bf g}}{|{\bf g}|}\cdot{\bf v}
+Fr\left(\rho^r-1\right)\int_{\Omega_n^s}\frac{{\bf g}}{|{\bf g}|}\cdot{\bf v} 
+\int_{\Gamma_N}\tilde{\bf h}\cdot{\bf v}
-\boxed{\mathcal{F}_{n+1/3}^{FSI}}.
\end{split}
\end{equation}
\end{itemize}
As noted above, expression (\ref{diffusion_explicit}) corresponds to a formulation of the classical explicit IFEM. The implicit IFEM may be expressed in a similar form, based upon the splitting scheme, but using the current value ${\bf u}_{n+2/3}$ to construct $\mathcal{F}_{n+2/3}^{FSI}$.

\begin{itemize}
\item Implicit forcing:
\begin{equation}\label{diffusion_implicit}
\begin{split}
&\int_{\Omega}\frac{{\bf u}_{n+2/3}-{\bf u}_{n+1/3}}{\Delta t} \cdot{\bf v}
+\frac{1}{2Re}\int_{\Omega}{\rm D}_n{\bf u}_{n+2/3}:{\rm D}_n{\bf v}  \\
=& Fr\int_{\Omega}\frac{{\bf g}}{|{\bf g}|}\cdot{\bf v}
+Fr\left(\rho^r-1\right)\int_{\Omega_n^s}\frac{{\bf g}}{|{\bf g}|}\cdot{\bf v} 
+\int_{\Gamma_N}\tilde{\bf h}\cdot{\bf v}
-\boxed{\mathcal{F}_{n+2/3}^{FSI}}.
\end{split}
\end{equation}
\end{itemize}

%%%%%%%%%%%%%%%%%%%%%%%%%%%
\section{One-field FDM}
\label{one_field_fdm}
It can be seen that the solid velocity is hidden in the nonlinear term ${\bf F}{\bf F}^T-{\bf I}$ in equation (\ref{weak_form1_time}) or (\ref{diffusion}). The difference between the one-field FDM and the explicit IFEM is how to treat this nonlinear term: the former extracts this hidden velocity by linearizing ${\bf F}{\bf F}^T-{\bf I}$ in the current configuration, while the latter evaluates this term as a force term on the right-hand side of the equation. In this section, we shall demonstrate how the nonlinear term ${\bf F}{\bf F}^T-{\bf I}$ is linearized in the one-field FDM, and expressed in terms of velocity in the current configuration. Also notice that this idea is not limited to the operator splitting. The splitting is just a specific implementation that allows us to express the IFEM and the one-field FDM in a similar form so as to facilitate direct comparison with each other.

Denoting ${\bf F}_t{\bf F}_t^T-{\bf I}$ by ${\bf s}_t$, and according to the definition of {\bf F} (\ref{definition_f}), ${{\bf s}_t}$ may be computed as follows:
\begin{equation}
{{\bf s}_t}
={\bf F}_t{\bf F}_t^T-{\bf I}
=\left(\nabla_{\bf X}{\bf x}_t\nabla_{\bf X}^T{\bf x}_t-{\bf I}\right).
\end{equation}
Using the chain rule, this last equation can also be expressed as:
\begin{equation}
{{\bf s}_t}
=\nabla_n{\bf x}_t \nabla_{\bf X}{\bf x}_n  \nabla_{\bf X}^T{\bf x}_n  \nabla_n^T{\bf x}_t
-{\bf I}
+\nabla_n{\bf x}_t   \nabla_n^T{\bf x}_t
-\nabla_n{\bf x}_t   \nabla_n^T{\bf x}_t
\end{equation}
or
\begin{equation}
{{\bf s}_t}
=\nabla_n{\bf x}_t   \nabla_n^T{\bf x}_t - {\bf I}
+\nabla_n{\bf x}_t\left(\nabla_{\bf X}{\bf x}_n  \nabla_{\bf X}^T{\bf x}_n-{\bf I}\right)\nabla_n^T{\bf x}_t.
\end{equation}
Then ${{\bf s}_t}$ can be expressed based on the previous coordinate ${\bf x}_n$ as follows:
\begin{equation}
{{\bf s}_t}
=\nabla_n{\bf x}_t   \nabla_n^T{\bf x}_t - {\bf I}
+\nabla_n{\bf x}_t {{\bf s}_{n}} \nabla_n^T{\bf x}_t.
\end{equation}
Using ${\bf x}_t={\bf x}_n+\Delta t{\bf u}_t$ (see Remark \ref{fraction_f} (\ref{fraction_f_eq})), the last equation can finally be expressed as:
\begin{equation}\label{solidstress}
\begin{split}
{{\bf s}_{t}}
&=\Delta t\left(\nabla_n{\bf u}_{t}+\nabla_n^T{\bf u}_{t}+\Delta t\nabla_n{\bf u}_{t}\nabla_n^T{\bf u}_{t}\right)
+{{\bf s}_n}   \\
&+\Delta t^2\nabla_n{\bf u}_{t} {{\bf s}_n} \nabla_n^T{\bf u}_{t}
+\Delta t\nabla_n{\bf u}_{t} {{\bf s}_n}
+\Delta t {{\bf s}_n} \nabla_n^T{\bf u}_{t}.
\end{split}
\end{equation}
There are two nonlinear terms in the last equation, which may be linearized as
\begin{equation}\label{linearization1}
\nabla_n{\bf u}_{t}\nabla_n^T{\bf u}_{t}
=\nabla_n{\bf u}_{t}\nabla_n^T{\bf u}_n
+\nabla_n{\bf u}_n\nabla_n^T{\bf u}_{t}
-\nabla_n{\bf u}_n\nabla_n^T{\bf u}_n,
\end{equation}
and
\begin{equation}\label{linearization2}
\nabla_n{\bf u}_{t} {{\bf s}_n} \nabla_n^T{\bf u}_{t}
=\nabla_n{\bf u}_{t} {{\bf s}_n} \nabla_n^T{\bf u}_n
+\nabla_n{\bf u}_n {{\bf s}_n} \nabla_n^T{\bf u}_{t}
-\nabla_n{\bf u}_n {{\bf s}_n} \nabla_n^T{\bf u}_n.
\end{equation}

Substituting ${\bf s}_{n+2/3}={{\bf F}_{n+2/3}{\bf F}_{n+2/3}^T-{\bf I}}$, using expression (\ref{solidstress}), (\ref{linearization1}) and (\ref{linearization2}), into diffusion step (\ref{diffusion}), we finally get the one-field FDM formulation as follows:
\begin{equation}\label{diffusion_fdm}
\begin{split}
&\int_{\Omega}\frac{{\bf u}_{n+2/3}-{\bf u}_{n+1/3}}{\Delta t} \cdot{\bf v}
+\frac{1}{2Re}\int_{\Omega}{\rm D}{\bf u}_{n+2/3}:{\rm D}{\bf v}  \\
+&\left(\rho^{r}-1\right)\int_{\Omega_n^s}\frac{{\bf u}_{n+2/3}-{\bf u}_n}{\Delta t} \cdot{\bf v}
+\frac{\mu^r-1}{2Re}\int_{\Omega_n^s}{\rm D}_n{\bf u}_{n+2/3}:{\rm D}_n{\bf v} \\
+& \boxed{\frac{\Delta t\tilde{c}_1}{2}\int_{\Omega_n^s}{\rm D}_n{\bf u}_{n+2/3}:{\rm D}_n{\bf v}} 
+ \boxed{\Delta t\tilde{c}_1\int_{\Omega_n^s}{\rm D}_n^1{\bf u}_{n+2/3}:\nabla_n{\bf v}}  \\
+& \boxed{\Delta t^2\tilde{c}_1\int_{\Omega_n^s}\left({\rm D}_n^2+{\rm D}_n^3\right){\bf u}_{n+2/3}:\nabla_n{\bf v}} \\
=& Fr\int_{\Omega}\frac{{\bf g}}{|{\bf g}|}\cdot{\bf v}
+Fr\left(\rho^r-1\right)\int_{\Omega_n^s}\frac{{\bf g}}{|{\bf g}|}\cdot{\bf v} 
+\int_{\Gamma_N}\tilde{\bf h}\cdot{\bf v}\\
-& \tilde{c}_1\int_{\Omega_n^s}{\bf s}_n:\nabla_n{\bf v} 
+ \boxed{\Delta t^2\tilde{c}_1\int_{\Omega_n^s}\left(\nabla_n{\bf u}_n\nabla_n^T{\bf u}_n\right):\nabla_n{\bf v}} \\
+ & \boxed{\Delta t^2\tilde{c}_1\int_{\Omega_n^s}\left(\nabla_n{\bf u}_n {\bf s}_n \nabla_n^T{\bf u}_n\right):\nabla_n{\bf v}}.
\end{split}
\end{equation}
In the above, the linear operators ${\rm D}_n^1$, ${\rm D}_n^2$ and ${\rm D}_n^3$ are defined as:
\begin{equation}
{\rm D}_n^1{\bf u}=\nabla_n{\bf u}{\bf s}_n+{\bf s}_n\nabla_n^T{\bf u},
\end{equation}
\begin{equation}
{\rm D}_n^2{\bf u}=\nabla_n{\bf u}\nabla_n^T{\bf u}_n+\nabla_n{\bf u}_n\nabla_n^T{\bf u},
\end{equation}
and
\begin{equation}
{\rm D}_n^3{\bf u}=\nabla_n{\bf u}{\bf s}_n\nabla_n^T{\bf u}_n+\nabla_n{\bf u}_n{\bf s}_n\nabla_n^T{\bf u}.
\end{equation}

\begin{remark}\label{not_trivial_terms}
When $\rho^r=\mu^r=1$, comparing equations (\ref{diffusion_explicit}) and (\ref{diffusion_fdm}), we see that the only difference between the one-field FDM and explicit IFEM is the additional boxed terms in equation (\ref{diffusion_fdm}) of $O(\Delta t)$ or $O(\Delta t^2)$ respectively. We shall demonstrate, with numerical tests, that these terms are not trivial at all: in fact they have a very positive stabilizing effect, such that a significantly larger time step may be adopted.
\end{remark}

%%%%%%%%%%%%%%%%%%%%%%%%%%%%%%%%%%%%%%%%%%
\section{Numerical experiments}
\label{sec_numerical_exs}
Having analyzed the proposed one-field FDM, and illustrated its close relationship with the IFEM family of methods, in this section we focus on validation of the three advantages, as claimed in Section \ref{sec_introduction}, of the one-field FDM compared with IFEM. We shall use the Least-squares method to approximate the convection step \cite{Wang_2017}. For the diffusion step in which FSI interaction happens, we use the finite element isoparametric interpolation functions to transfer information between the solid mesh and fluid mesh for both the one-field FDM and IFEM. The finite element interpolation function is suggested to be capable of producing sharper interfaces than the traditional discretized Dirac delta function or the reproducing kernel function in \cite{Wang_2009}. The pressure step is a ``degenerate'' Stokes equation, and we shall use the standard $Taylor$-$Hood$ element to discretize this step. Gravity will not be considered in this paper, so the Froude number (\ref{froude}) will be zero ($Fr=0$) in each of the following tests.

\subsection{Lid-driven cavity flow with a solid disc}
\begin{figure}[h!]
\centering
\includegraphics[width=2.5in,angle=0]{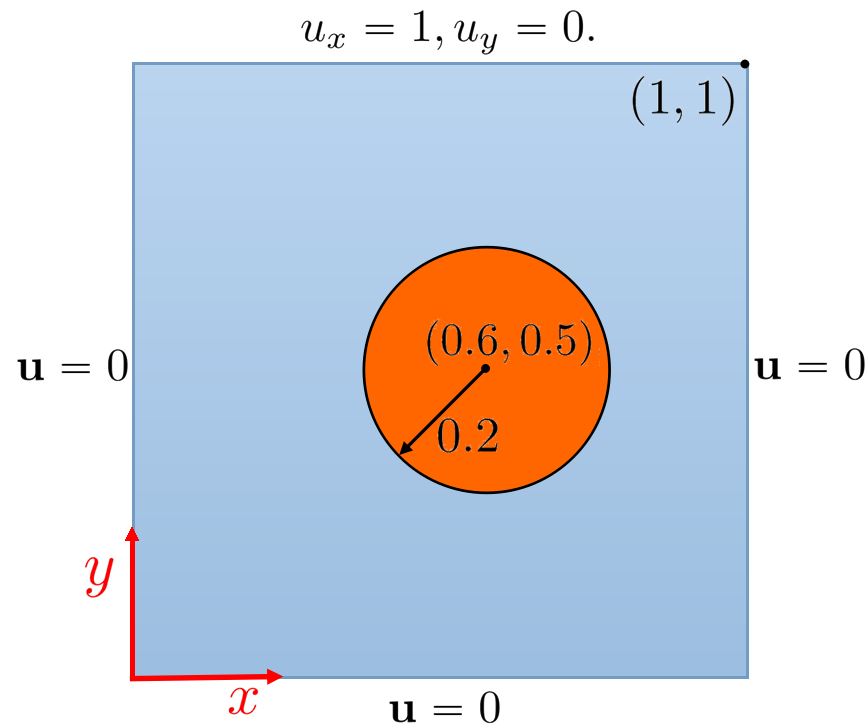}
\captionsetup{justification=centering}
\caption {\scriptsize Computational domain and boundary conditions for the test problem of a lid-driven cavity flow with a solid disc.} 
\label{Sketch_cavity}		
\end{figure}

\begin{figure}[h!]
	\begin{minipage}[t]{0.5\linewidth}
		\centering
		\includegraphics[width=2.2in,angle=0]{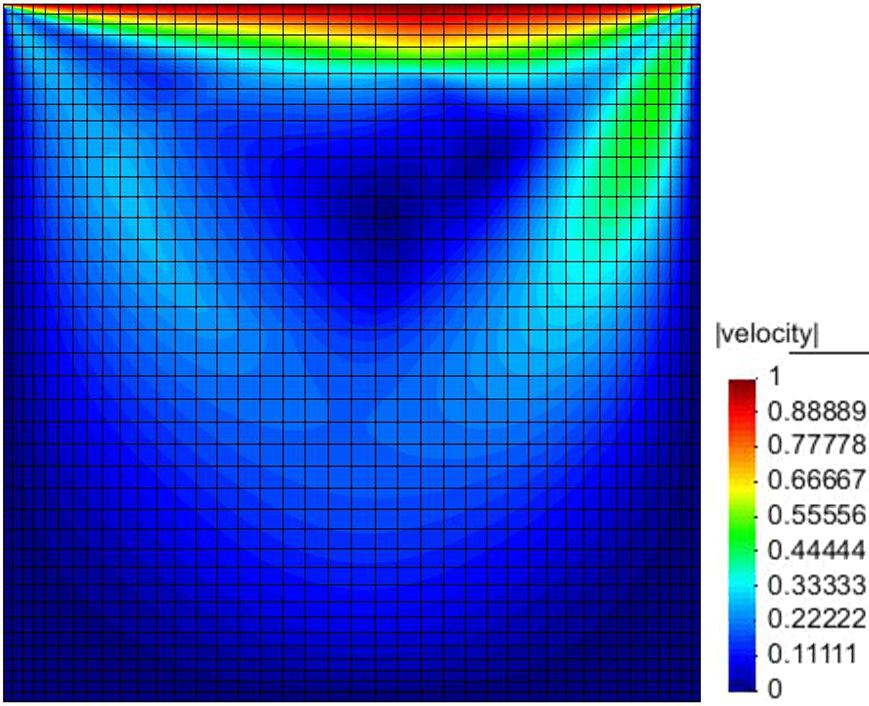}
		\captionsetup{justification=centering}
		\caption* {\scriptsize (a) Velocity on fluid mesh.} 
	\end{minipage}
	\begin{minipage}[t]{0.5\linewidth}
		\centering
		\includegraphics[width=2.0in,angle=0]{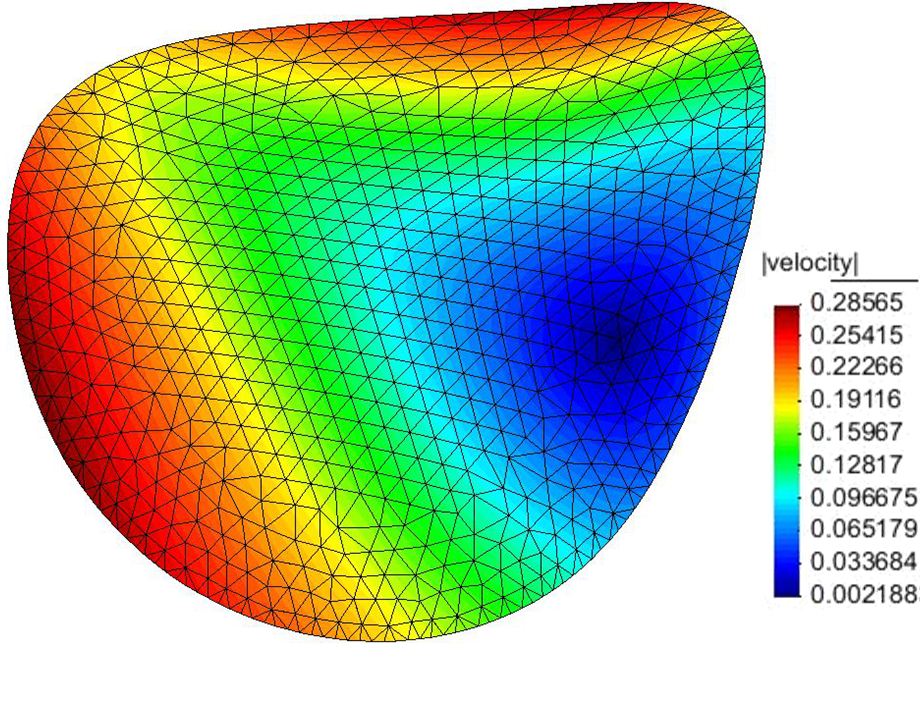}
		\captionsetup{justification=centering}
		\caption* {\scriptsize (b) Velocity on solid mesh.} 
	\end{minipage} 		
	\captionsetup{justification=centering}
	\caption {\scriptsize  Distribution of the velocity norm for Parameter set 1 at $t=10$. The position of the solid in the cavity can be seen from Figure \ref{cavity_solid_deformation_mu1} (d). } 
	\label{cavity_v_fluid_solid}
\end{figure}

This example is taken from papers \cite{Zhao_2008,Wang_2009}, in which IFEM approaches are adopted. A sketch of the problem and boundary conditions are shown in Figure \ref{Sketch_cavity}. We consider the parameter sets displayed in Table \ref{Parameter_sets_cavity}, and all the simulations use a sufficiently small time step to ensure stability: parameter set 1 (very soft solid) uses $\Delta t=1.0\times 10^{-3}$ and all other tests have a time step of $\Delta t=5.0\times 10^{-3}$. To illustrate the meshes that we use, velocity norms on the fluid mesh ($40\times 40$ biquadratic quadrilaterals) and solid mesh (1373 linear triangles with 771 nodes) for Parameter set 1 are presented in Figure \ref{cavity_v_fluid_solid}.

\begin{table}[h!]
\centering
\begin{tabular}{|c|c|c|c|c|c|}
\hline
Parameter sets  & $Re$ & $\tilde{c}_1$ & $\rho^r$ & $\mu^r$ \\
\hline 
Parameter 1 & $100$ & $0.1$ & $1$ & $1$  \\ 
\hline 
Parameter 2 & $100$ & $1$ & $1$ & $1$ \\ 
\hline 
Parameter 3 & $100$ & $1$ & $2$ & $1$ \\ 
\hline 	
Parameter 4 & $100$ & $1$ & $0.5$ & $1$ \\ 
\hline
Parameter 5 & $100$ & $1$ & $1$ & $1.5$ \\ 
\hline 	
Parameter 6 & $100$ & $1$ & $1$ & $2$  \\ 
\hline 	
Parameter 7 & $500$ & $0.5$ & $2$ & $2$ \\ 
\hline 
\end{tabular} \\                                       	   
\captionsetup{justification=centering}
\caption{Parameter sets for lid-driven cavity flow with a solid disc (the first parameter set is used in \cite{Zhao_2008,Wang_2009}).}
\label{Parameter_sets_cavity}
\end{table}

For all these tests, we aim to run to $t=10$. However for the Parameter sets 3 to 5 ($\rho^r\ne 1$), our implementation of the implicit IFEM scheme cannot reach $t=10$. Therefore, in these cases, we show a comparison shortly before the IFEM breaks down. The following criterion is used to check whether the implicit IFEM converges.
\begin{equation}\label{stop_criterion}
error=\frac{\|{\bf u}_{k+1}-{\bf u}_{k}\|_{\Omega_n^s}}{\|{\bf u}_{k}\|_{\Omega_n^s}} < tol,
\end{equation}
where ${\bf u}_k$ and ${\bf u}_{k+1}$ are the iterative values of the last and current step respectively, and $tol=10^{-6}$ is the convergence tolerance used in our tests.

\begin{figure}[h!]
	\begin{minipage}[t]{0.5\linewidth}
		\centering
		\includegraphics[width=2.2in,angle=0]{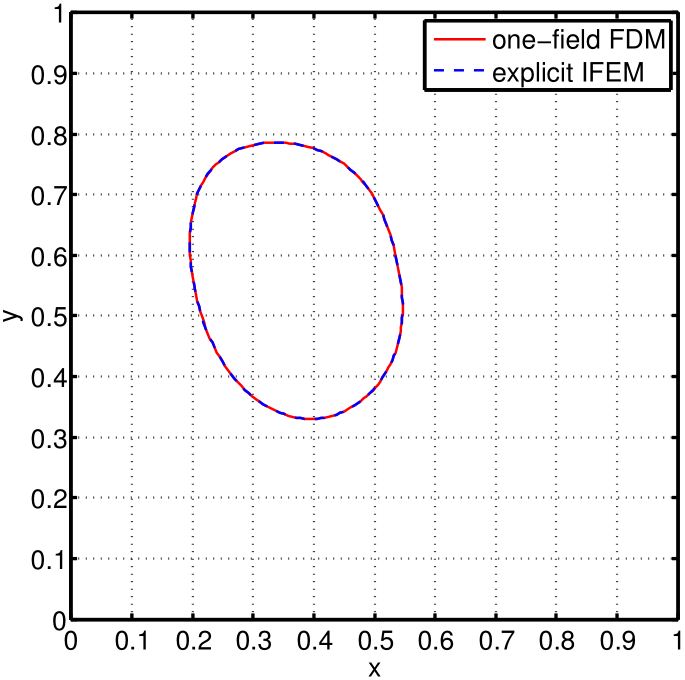}
		\captionsetup{justification=centering}
		\caption* {\scriptsize (a) $t=2.4$,} 
	\end{minipage}
	\begin{minipage}[t]{0.5\linewidth}
		\centering
		\includegraphics[width=2.2in,angle=0]{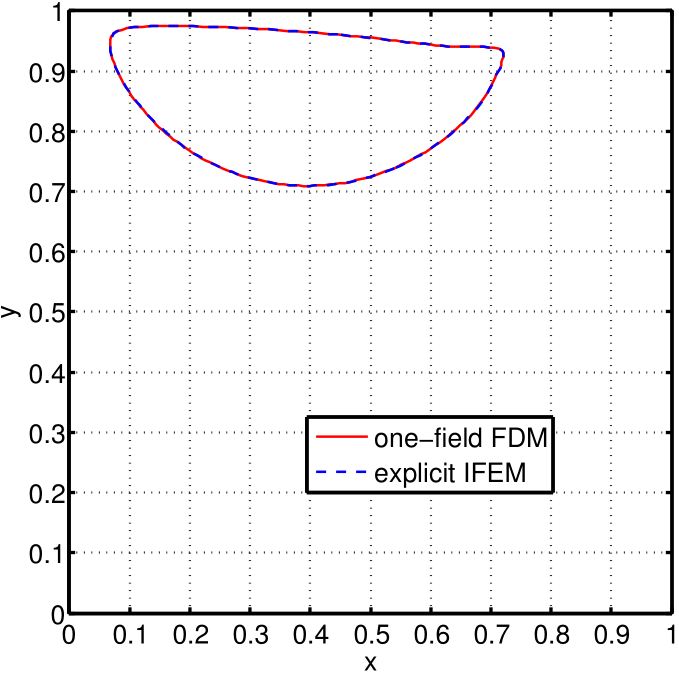}
		\captionsetup{justification=centering}
		\caption* {\scriptsize (b) $t=4.7$,} 
	\end{minipage} 
	\begin{minipage}[t]{0.5\linewidth}
		\centering
		\includegraphics[width=2.2in,angle=0]{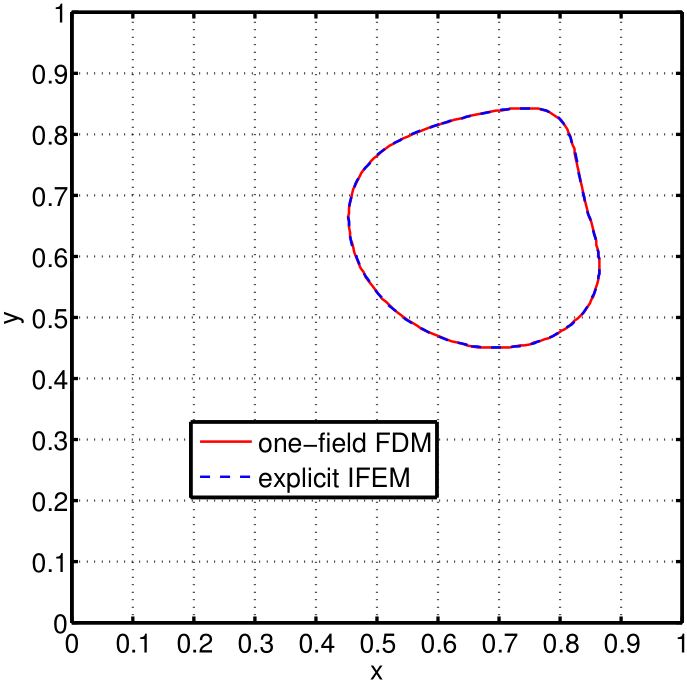}
		\captionsetup{justification=centering}
		\caption* {\scriptsize (c) $t=7.0$,} 
	\end{minipage}
	\begin{minipage}[t]{0.5\linewidth}
		\centering
		\includegraphics[width=2.2in,angle=0]{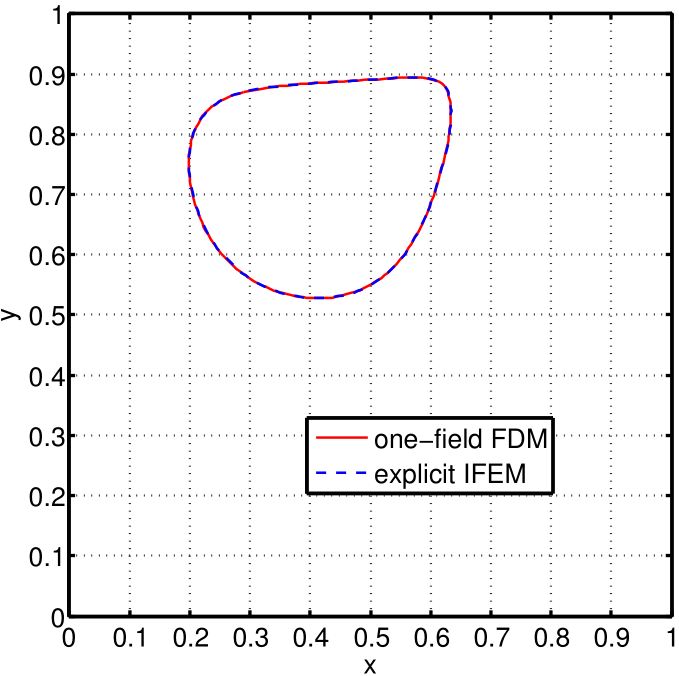}
		\captionsetup{justification=centering}
		\caption* {\scriptsize (d) $t=10$.} 
	\end{minipage}   		
	\captionsetup{justification=centering}
	\caption {\scriptsize  Solid deformation for Parameter set 1. These results are identical to Fig. 10 in \cite{Wang_2009}. The figures show that the one-field FDM and the explicit IFEM present the same results in the case of $\rho^r=\mu^r=1$. The $l^2$ norm of velocity vectors on the solid mesh at $t=10$: $\|{\bf u}_{\rm IFEM}\|=4.80955$, $\|{\bf u}_{\rm FDM}\|=4.80087$ and $\|{\bf u}_{\rm IFEM}-{\bf u}_{\rm FDM}\|=0.07399$.}
	\label{cavity_solid_deformation_mu1}
\end{figure}
\begin{figure}[h!]
\centering
\includegraphics[width=5.0 in,angle=0]{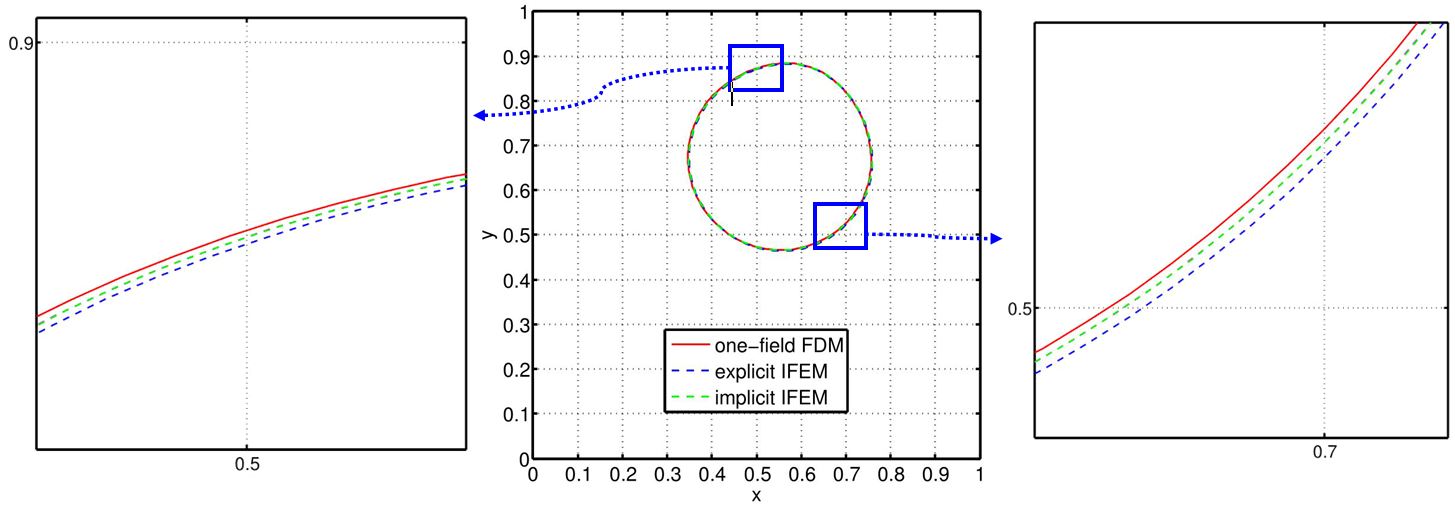}
\captionsetup{justification=centering}
\caption {\scriptsize Solid deformation for Parameter set 2 at $t=10$. The $l^2$ norm of displacement vectors on the solid mesh: $\|{\bf d}_{\rm exIFEM}-{\bf d}_{\rm FDM}\|=0.13423$ and $\|{\bf d}_{\rm imIFEM}-{\bf d}_{\rm FDM}\|=0.12248$}.
\label{cavity_solid_deformation_mu2}
\end{figure}

The first two parameter sets are simple cases because $\rho^r=\mu^r=1$. We can see from Figure \ref{cavity_solid_deformation_mu1} that the one-field FDM and the explicit IFEM present almost identical results in the case of a very soft solid ($\tilde{c}_1=0.1$), both of which are themselves indistinguishable from the published results in the literature \cite{Wang_2009}. Figure \ref{cavity_solid_deformation_mu2} shows the solid deformation for a slightly harder disc ($\tilde{c}_1=1$). Although the explicit IFEM, implicit IFEM and the one-field FDM present very similar results, a close look at the interface shape shows that the solution of the one-field FDM is almost identical to that of the implicit IFEM, and different from the explicit IFEM.

\begin{figure}[h!]
\centering
\includegraphics[width=2.5 in,angle=0]{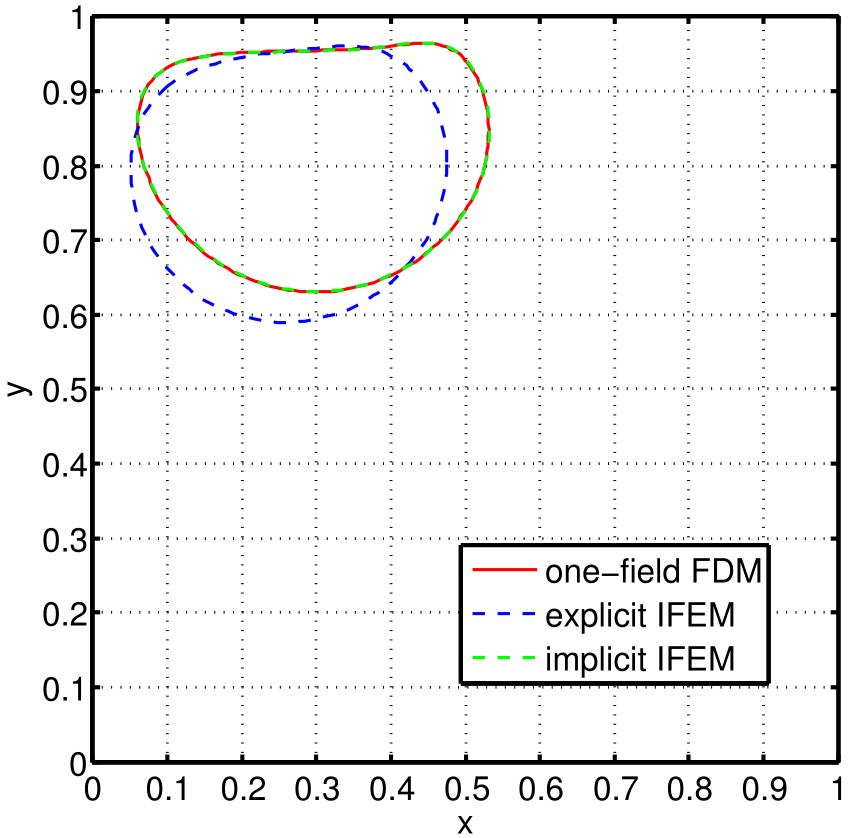}
\captionsetup{justification=centering}
\caption {\scriptsize Solid deformation for Parameter set 3 at $t=4.4$.} 
\label{cavity_solid_deformation_parameter3}
\end{figure}

\begin{figure}[h!]
\centering
\includegraphics[width=4.0 in,angle=0]{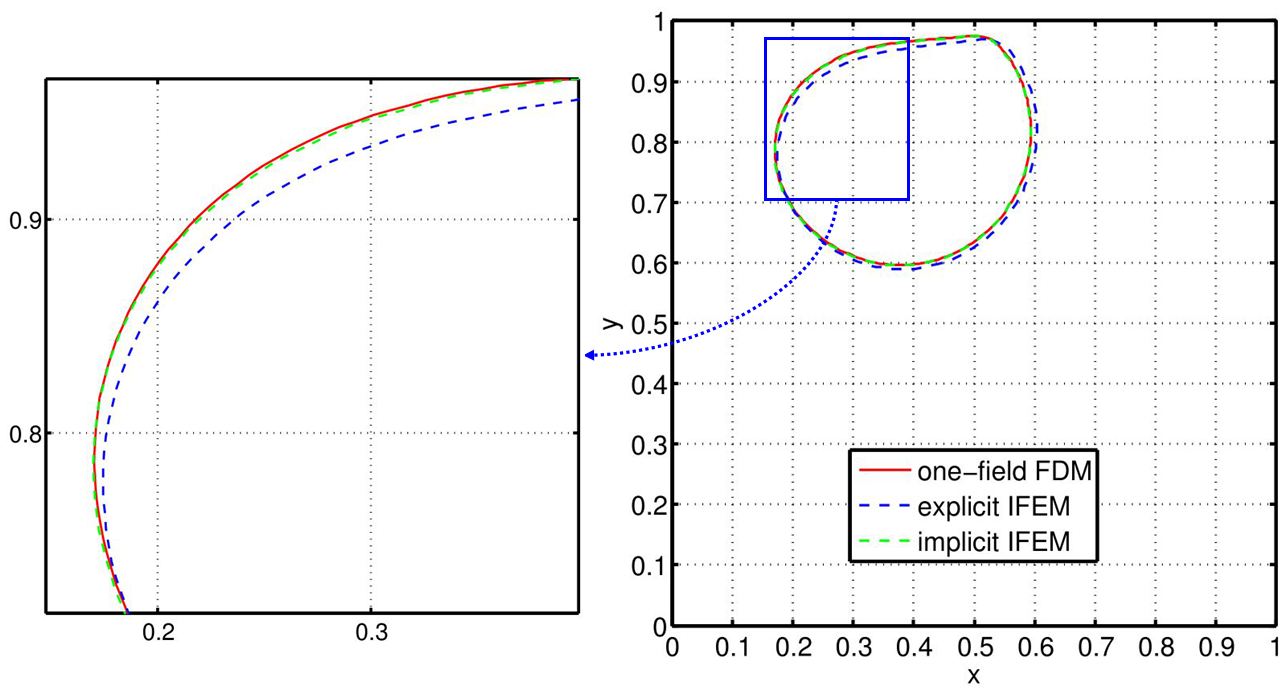}
\captionsetup{justification=centering}
\caption {\scriptsize Solid deformation for Parameter set 4 at $t=5.0$. The result of the one-field FDM is identical to the result of the implicit IFEM.}
\label{cavity_solid_deformation_parameter4}
\end{figure}

The case of a larger solid density ($\rho^r=2$) and a smaller solid density ($\rho^r=0.5$) are tested by Parameter sets 3 and 4 respectively. Both results (Figures \ref{cavity_solid_deformation_parameter3} and \ref{cavity_solid_deformation_parameter4}) show that the one-field FDM and the fully converged implicit IFEM have almost the same accuracy. Furthermore, in neither IFEM case do the results converge up to $t=10$ when using the same time step as the one-field FDM: $\Delta t=5.0\times 10^{-3}$. The explicit IFEM uses velocities from the previous two time steps to compute the effect of the solid: it can be seen from Figure \ref{cavity_solid_deformation_parameter3} that the disc moves more slowly using this explicit IFEM. We also note that reducing the time step cannot solve the problem in this case, because the temporal term in the FSI force $\mathcal{F}_t^{FSI}$ (\ref{force_fsi}) becomes huge and has a negative effect on the stability. Figure \ref{cavity_solid_deformation_parameter4} demonstrates similar problems for the explicit IFEM, but the disc using the explicit IFEM moves faster than the one-field FDM or the implicit IFEM.

\begin{figure}[h!]
\centering
\includegraphics[width=2.0 in,angle=0]{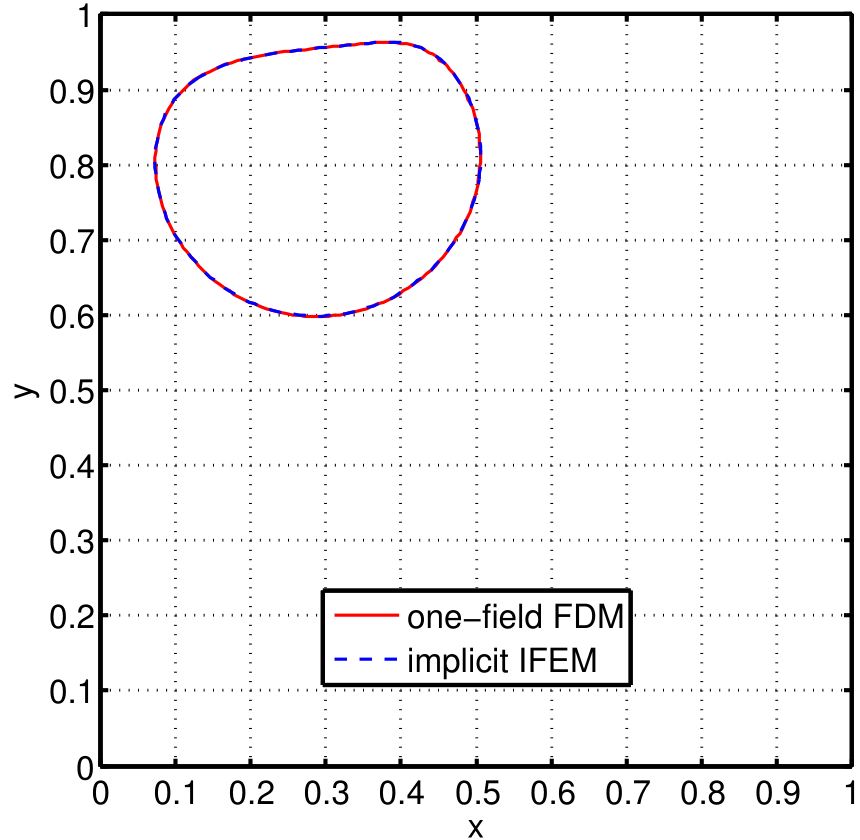}
\captionsetup{justification=centering}
\caption {\scriptsize Solid deformation for Parameter set 5 at $t=4.2$. The result of the one-field FDM is identical to the result of the implicit IFEM.}
\label{cavity_solid_deformation_parameter5}
\end{figure}

\begin{figure}[h!]
\begin{minipage}[t]{0.5\linewidth}
\centering
\includegraphics[width=2.4in,angle=0]{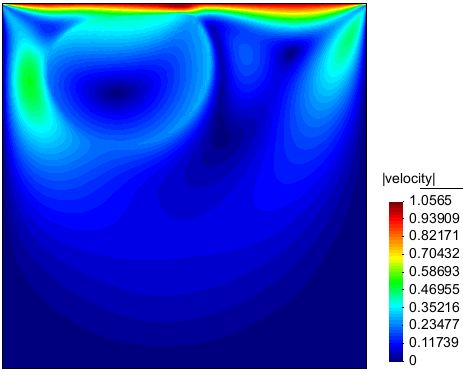}
\captionsetup{justification=centering}
\caption* {\scriptsize (a) Velocity in the background domain.} 
\end{minipage}
\begin{minipage}[t]{0.5\linewidth}
\centering
\includegraphics[width=2.0in,angle=0]{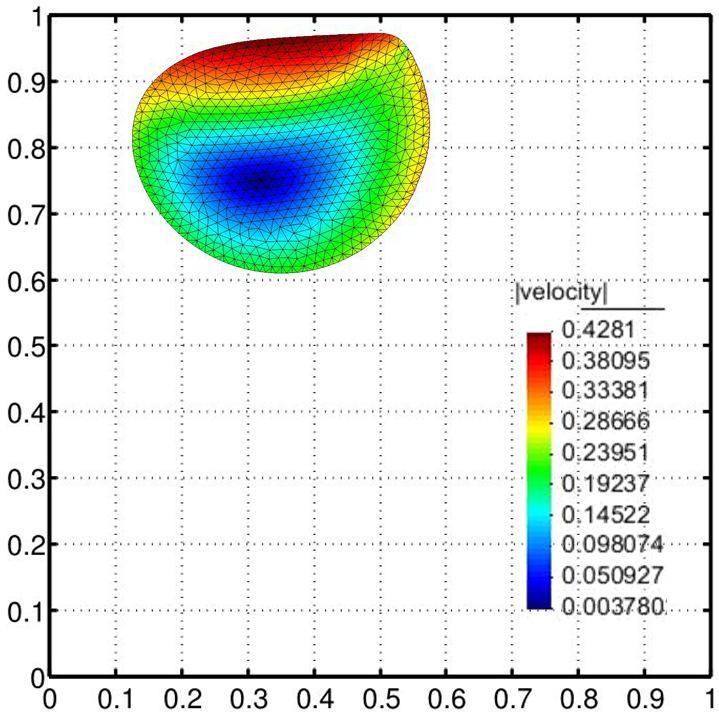}
\captionsetup{justification=centering}
\caption* {\scriptsize (b) Velocity on the solid mesh.} 
\end{minipage} 		
\captionsetup{justification=centering}
\caption {\scriptsize  Distribution of the velocity norm for Parameter set 6 at $t=5$, using the one-field FDM. The disc arrives at top of the cavity.} 
\label{cavity_v_t_5}
\end{figure}

\begin{figure}[h!]
\begin{minipage}[t]{0.5\linewidth}
\centering
\includegraphics[width=2.4in,angle=0]{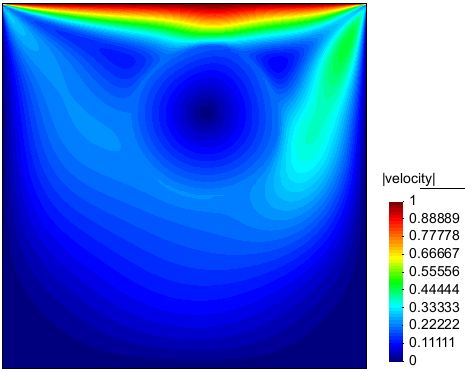}
\captionsetup{justification=centering}
\caption* {\scriptsize (a) Velocity in the background domain.} 
\end{minipage}
\begin{minipage}[t]{0.5\linewidth}
\centering
\includegraphics[width=2.0in,angle=0]{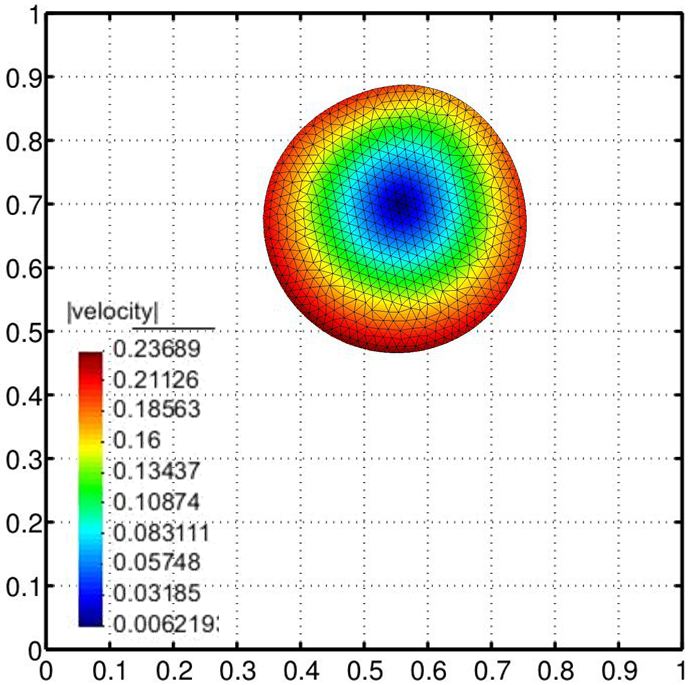}
\captionsetup{justification=centering}
\caption* {\scriptsize (b) Velocity on the solid mesh.} 
\end{minipage} 		
\captionsetup{justification=centering}
\caption {\scriptsize  Distribution of the velocity norm for Parameter set 6 at $t=10$, using the one-field FDM.} 
\label{cavity_v_t_10}
\end{figure}

Parameter sets 5 ($\mu^r=1.5$) and 6 ($\mu^r=2$) are included to show the case of different viscosities between fluid and solid. It can be seen from Figure \ref{cavity_solid_deformation_parameter5} that the results of the one-field FDM and the implicit IFEM match very well. Using the selected time step, our IFEM implementations break down after the first few time steps when testing Parameter set 6, therefore we only show the results of the one-field FDM in Figures \ref{cavity_v_t_5} and \ref{cavity_v_t_10} (as future test cases for others to validate against).

\begin{figure}[h!]
\begin{minipage}[t]{0.5\linewidth}
\centering
\includegraphics[width=2.4in,angle=0]{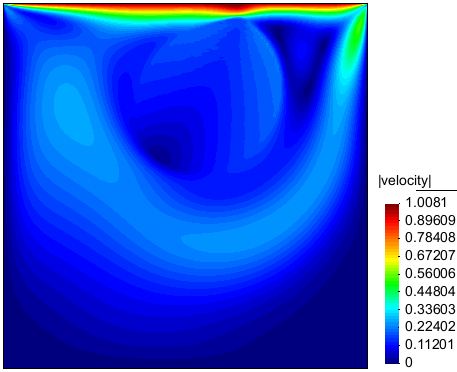}
\captionsetup{justification=centering}
\caption* {\scriptsize (a) Velocity in the background domain.} 
\end{minipage}
\begin{minipage}[t]{0.5\linewidth}
\centering
\includegraphics[width=2.0in,angle=0]{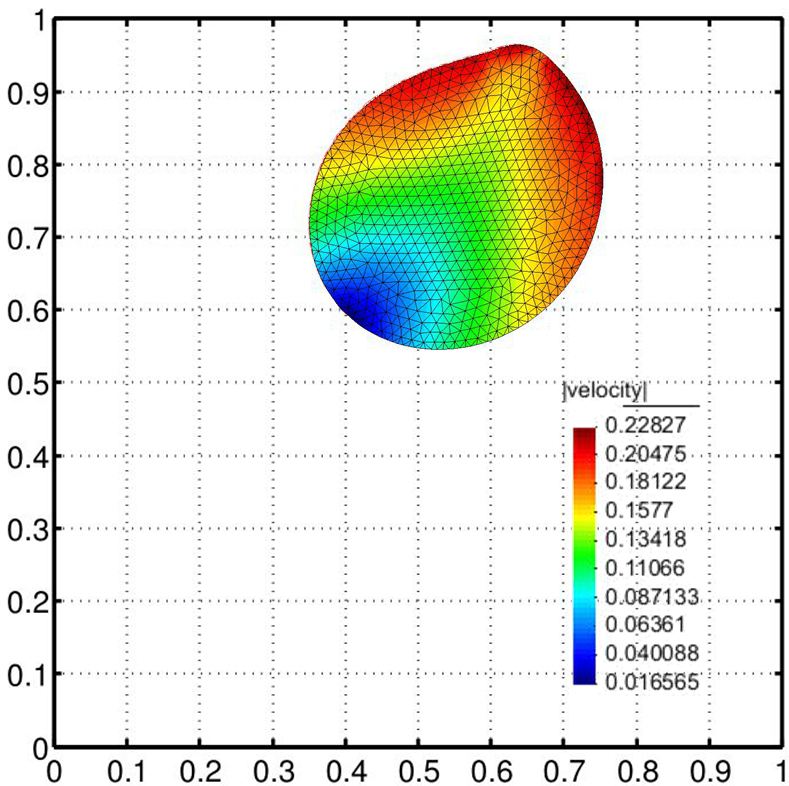}
\captionsetup{justification=centering}
\caption* {\scriptsize (b) Velocity on solid mesh.} 
\end{minipage} 		
\captionsetup{justification=centering}
\caption {\scriptsize  Distribution of the velocity norm for Parameter set 7 at $t=7.4$, using the one-field FDM. The disc arrives at top of the cavity.}  
\label{cavity_v_t_74}
\end{figure}

We purposely choose Parameter set 7 to be difficult, with large Reynolds number and differences in viscosity and density between fluid and solid. The one-field FDM can stably run up to $t=10$. We show the result in Figure \ref{cavity_v_t_74} when the solid disc arrives at the top of the cavity.

\begin{remark}
Notice that we have not considered the case of $\mu^r<1$, because we find that all the three methods (one-field FDM, explicit IFEM and implicit IFEM) may be unstable when $\Delta t\to 0$. However we shall not discuss this stability issue in more detail here as it is not the primary purpose of this paper. Please refer to \cite{Wang_2019,yongxing2018} for stability analysis.
\end{remark}

%%%%%%%%%%%%%%%%%%
\subsection{Oscillating leaflet in a channel}
\label{subsec_leaflet}
\begin{figure}[h!]
\centering
\includegraphics[width=5in,angle=0]{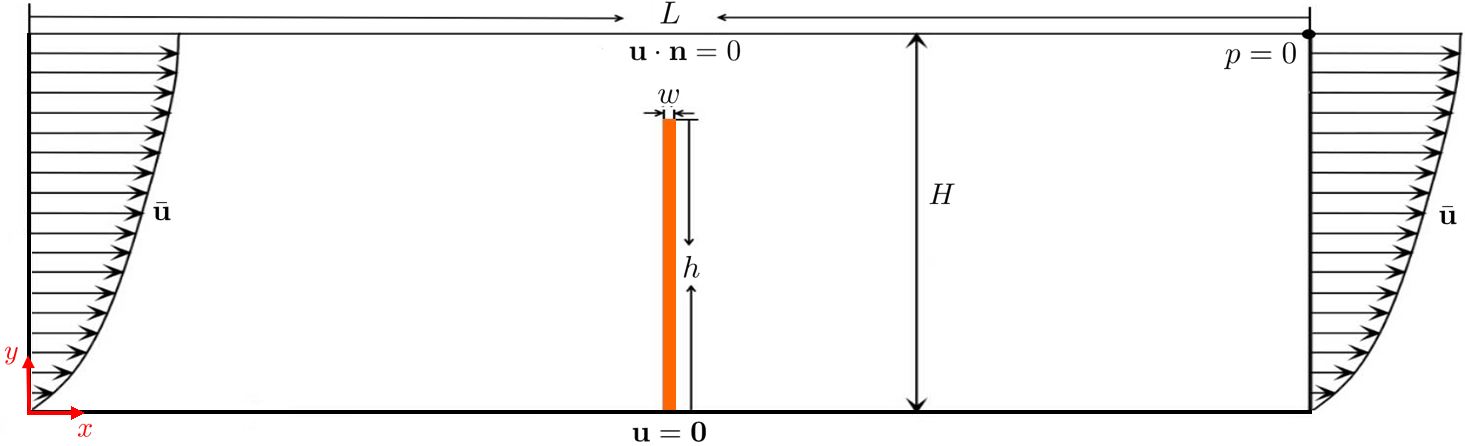}
\captionsetup{justification=centering}
\caption {\scriptsize Computational domain and boundary conditions for the oscillating leaflet.} 
\label{leaflet_geometry}
\end{figure}

This numerical example is taken from \cite{Yu_2005,baaijens2001fictitious,Kadapa_2016}. The computational domain is a channel of dimension $L\times H$, with a leaflet of dimension $h\times w$ initially located across the channel at its midpoint as shown in Figure \ref{leaflet_geometry}. A periodic flow condition is prescribed on the inlet and outlet boundaries, given by
\begin{equation}
\bar{u}_x=1.5y\left(2-y\right)sin\left(2\pi t/T\right),\quad \bar{u}_y=0,
\end{equation}
with $T$ being the dimensionless period of the input flow and equal to 10. In this test, $L=4$, $H=1$, $h=0.8$ and $w=0.0212$.

\begin{figure}[h!]
\begin{minipage}[t]{0.5\linewidth}
\centering
\includegraphics[width=2.0in,angle=0]{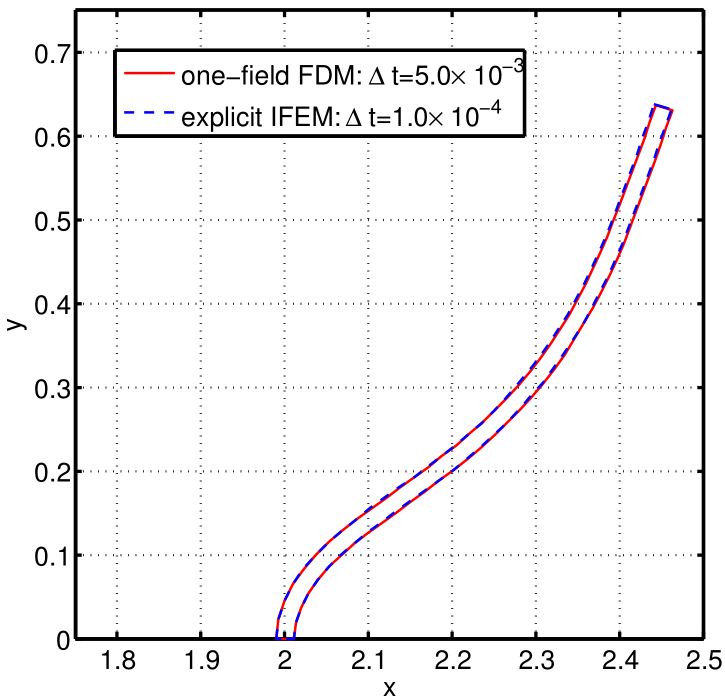}
\captionsetup{justification=centering}
\caption* {\scriptsize (a) $t/T=0.1$, $\|{\bf u}\|_r=8.62548\times 10^{-3}$}. 
\end{minipage}
\begin{minipage}[t]{0.5\linewidth}
\centering
\includegraphics[width=2.3in,angle=0]{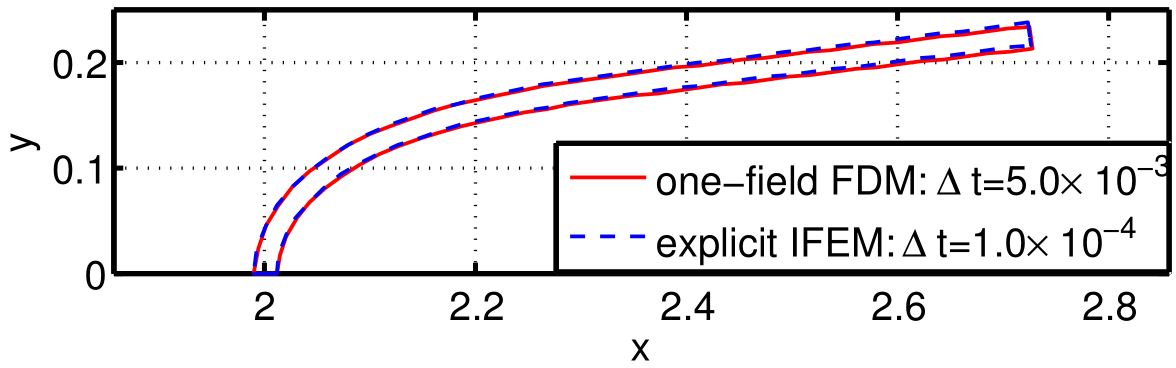}
\captionsetup{justification=centering}
\caption* {\scriptsize (b) $t/T=0.2$, $\|{\bf u}\|_r=4.52988\times 10^{-2}$}. 
\end{minipage} 
\begin{minipage}[t]{0.5\linewidth}
\centering
\includegraphics[width=2.2in,angle=0]{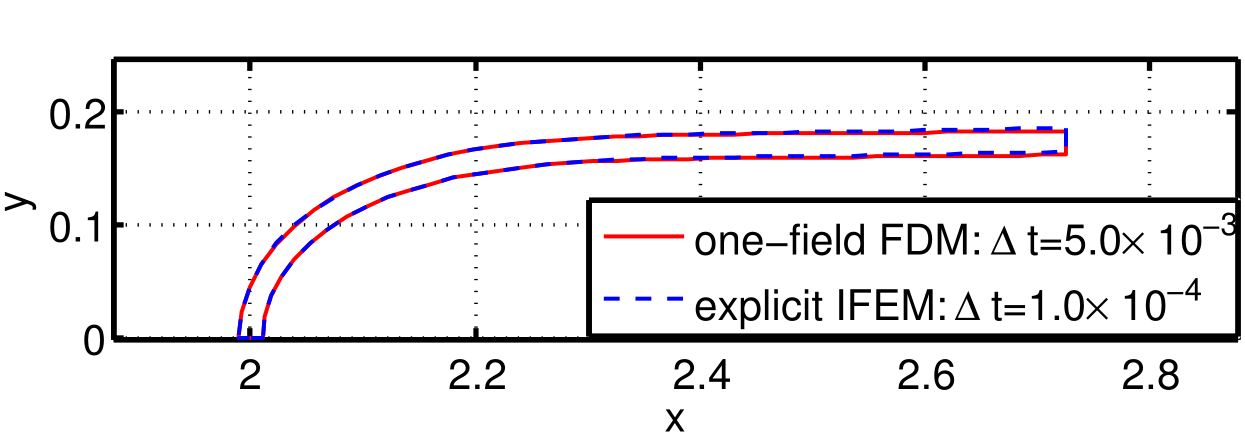}
\captionsetup{justification=centering}
\caption* {\scriptsize (c) $t/T=0.3$, $\|{\bf u}\|_r=4.12284\times 10^{-2}$}
\end{minipage}
\begin{minipage}[t]{0.5\linewidth}
\centering
\includegraphics[width=2.3in,angle=0]{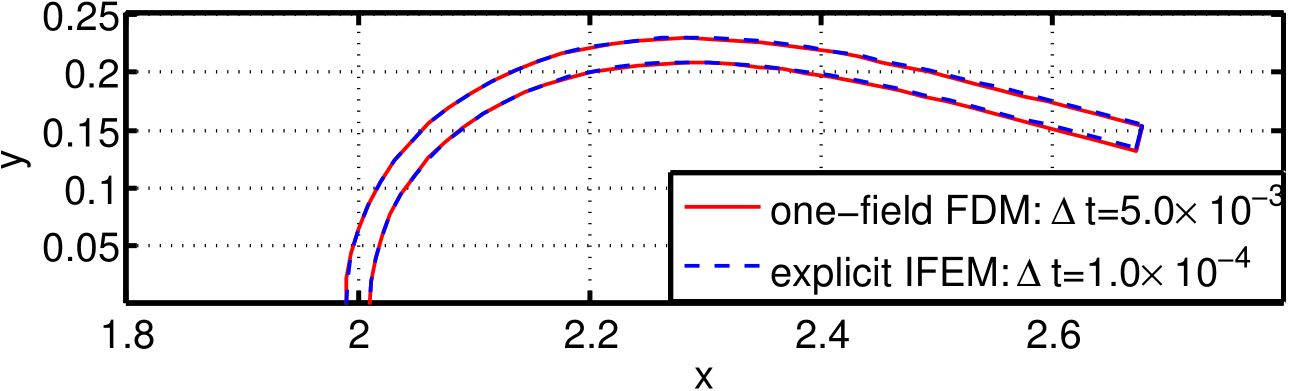}
\captionsetup{justification=centering}
\caption* {\scriptsize (d) $t/T=0.4$, $\|{\bf u}\|_r=1.03842\times 10^{-2}$}.
\end{minipage} 	
\captionsetup{justification=centering}
\caption {\scriptsize Leaflet deformation at different times using Parameter set 1. Comparison between the one-field FDM and the explicit IFEM shows excellent agreement with \cite{Yu_2005}. Error measured by the $l^2$ norm of velocity difference: $\|{\bf u}\|_r=\|{\bf u}_{\rm FDM}-{\bf u}_{\rm IFEM}\| / \|{\bf u}_{\rm FDM}\|$.}
\label{leaflet_defomration_parameter1}
\end{figure}

The leaflet is approximated with 152 linear triangles with 116 nodes, and the fluid mesh is made up of $189\times47$ biquadratic squares with 36005 nodes. We extend parameters ($\rho^r=1$) used in the above three publications to two more general cases as shown in Table \ref{Parameter_sets_leaflet}. Using the first group of parameters, we demonstrate that the one-field FDM can use a time step of $\Delta t=5.0\times 10^{-3}$ while the explicit IFEM has to use a time step of $\Delta t=1.0\times 10^{-4}$ in order to remain stable (we reduce the time by a half over consecutive tests to check the convergence until finding a stable time step $\Delta t=7.8125\times 10^{-5}$, and then slightly increasing it we find $\Delta t=1.0\times 10^{-4}$ is also stable). However both simulations lead to almost identical results, as shown in Figure \ref{leaflet_defomration_parameter1}, which match the results of Fig. 3 in \cite{Yu_2005}. The reason for the difference in $\Delta t$ is due to the additional stabilizing terms added in the one-field FDM, as discussed in Remark \ref{not_trivial_terms}. Also notice that we have to use the same time step for both our implicit and explicit IFEM schemes in order to converge for this example, because of the huge forcing term on the right-hand side in equation (\ref{diffusion_explicit}) and (\ref{diffusion_implicit}).
\begin{table}[h!]
	\centering
	\begin{tabular}{|c|c|c|c|c|}
		\hline
		Parameter sets  & $Re$ & $\tilde{c}_1$ & $\rho^r$ & $\mu^r$ \\
		\hline 
		Parameter 1 & $100$ & $1000$ & $1$ & $1$  \\ 
		\hline 
		Parameter 2 & $100$ & $1000$ & $1.2$ & $1$  \\ 
		\hline 
		Parameter 3 & $300$ & $10000$ & $2$ & $2$  \\ 
		\hline 
	\end{tabular} \\                                       	   
	\captionsetup{justification=centering}
	\caption{Parameter sets for the oscillating leaflet in a channel.}
	\label{Parameter_sets_leaflet}
\end{table}

\begin{figure}[h!]
	\begin{minipage}[t]{0.5\linewidth}
		\includegraphics[width=2.0in,angle=0]{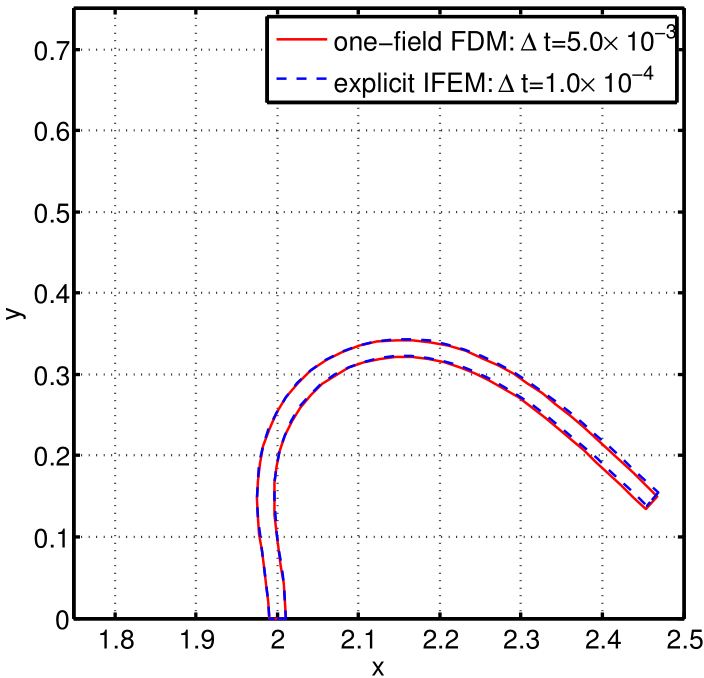}
		\captionsetup{justification=centering}
		\caption* {\scriptsize (e) $t/T=0.5$, $\|{\bf u}\|_r=9.28756\times 10^{-3}$}.
	\end{minipage}
	\begin{minipage}[t]{0.5\linewidth}
		\centering
		\includegraphics[width=2.1in,angle=0]{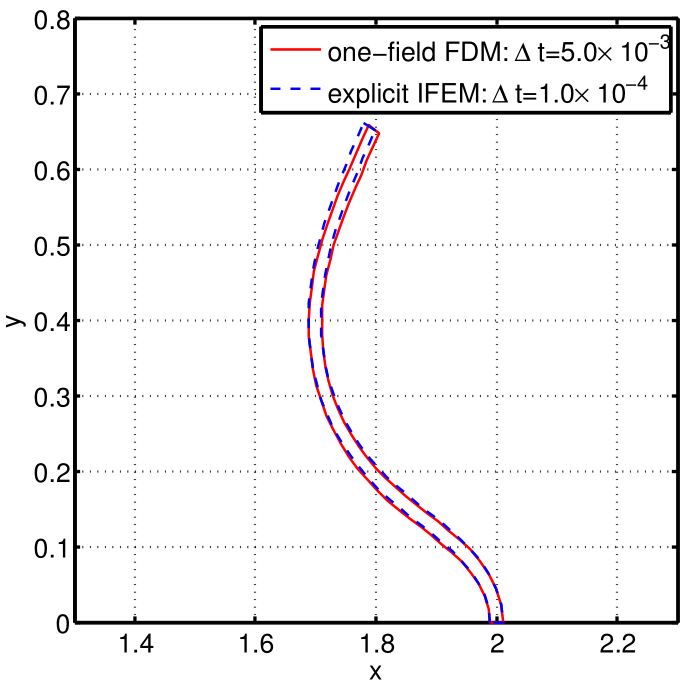}
		\captionsetup{justification=centering}
		\caption* {\scriptsize (f) $t/T=0.6$, $\|{\bf u}\|_r=3.03994\times 10^{-2}$}.
	\end{minipage} 
	\begin{minipage}[t]{0.5\linewidth}
		\centering
		\includegraphics[width=2.3in,angle=0]{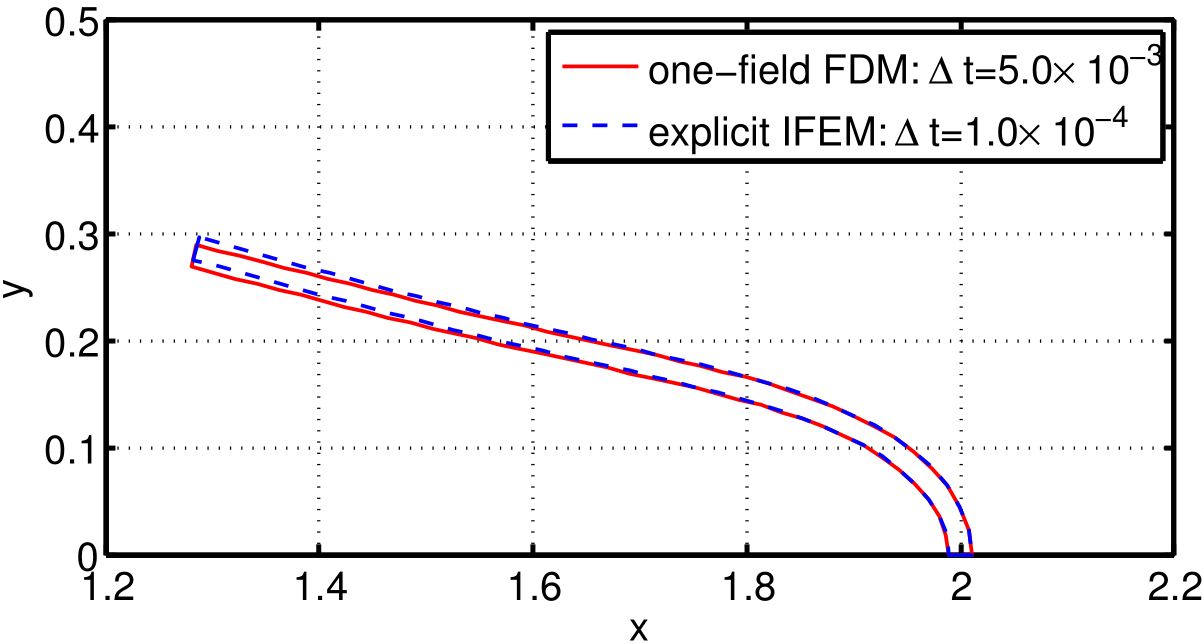}
		\captionsetup{justification=centering}
		\caption* {\scriptsize (g) $t/T=0.7$, $\|{\bf u}\|_r=5.85640\times 10^{-2}$}.
	\end{minipage}
	\begin{minipage}[t]{0.5\linewidth}
		\centering
		\includegraphics[width=2.3in,angle=0]{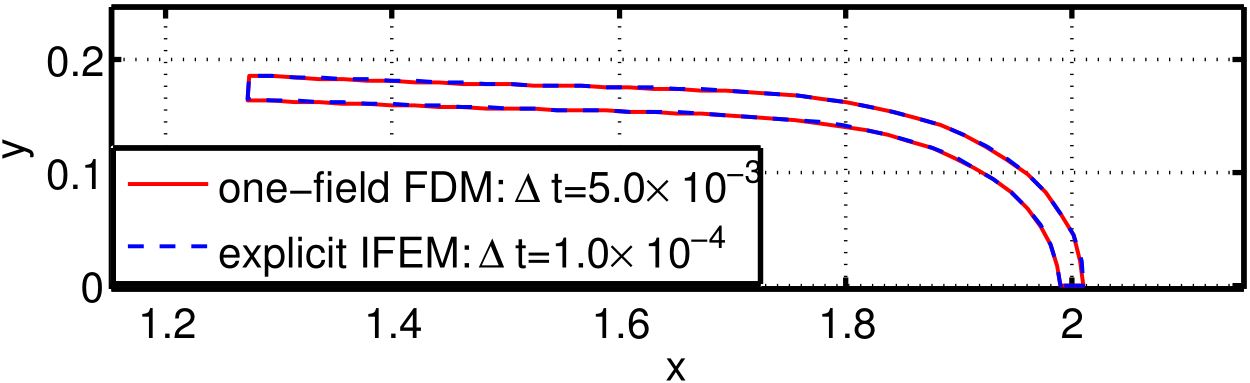}
		\captionsetup{justification=centering}
		\caption* {\scriptsize (h) $t/T=0.8$, $\|{\bf u}\|_r=2.54594\times 10^{-2}$}.
	\end{minipage} 
	\captionsetup{justification=centering}
	\caption* {Figure \ref{leaflet_defomration_parameter1} (continued).}  
\end{figure}

\begin{figure}[h!]
	\begin{minipage}[t]{0.5\linewidth}
		\centering
		\includegraphics[width=2.1in,angle=0]{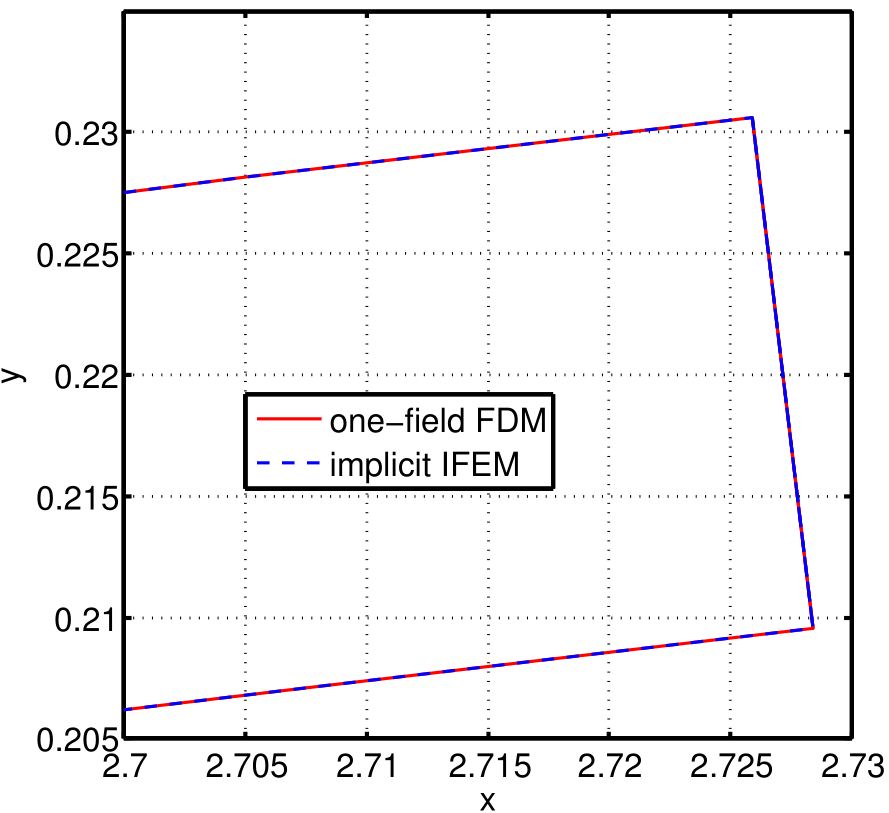}
		\captionsetup{justification=centering}
		\caption* {\scriptsize (a) $t/T=0.2$} 
	\end{minipage}
	\begin{minipage}[t]{0.5\linewidth}
		\centering
		\includegraphics[width=2.0in,angle=0]{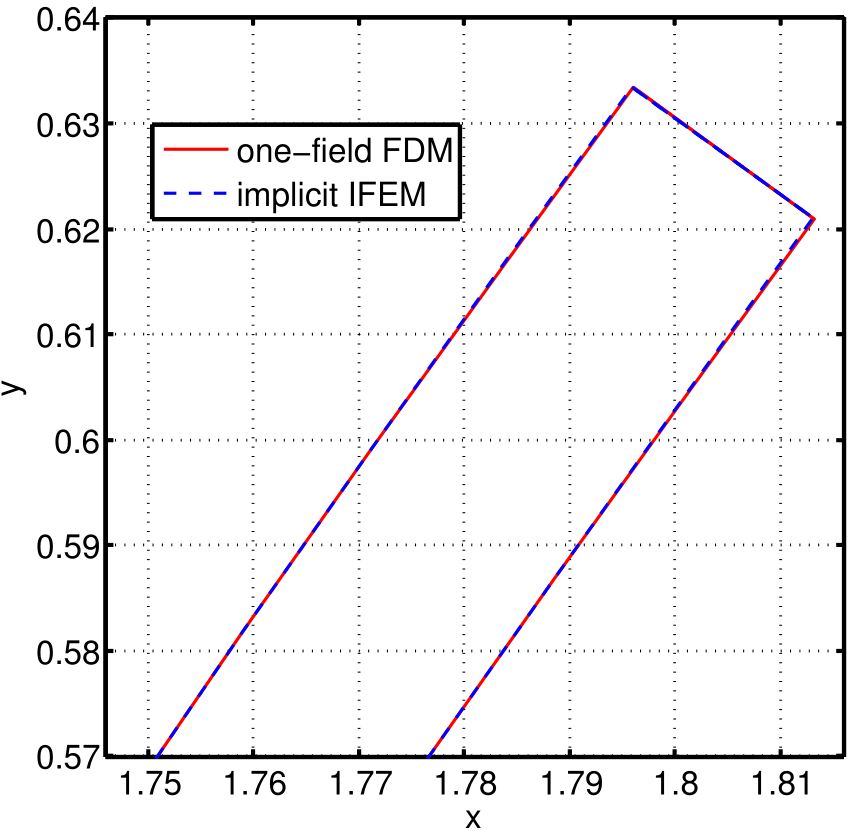}
		\captionsetup{justification=centering}
		\caption* {\scriptsize (b) $t/T=0.6$.} 
	\end{minipage} 
	\captionsetup{justification=centering}
	\caption {\scriptsize Deformation at the tip of the leaflet for Parameter set 2. Comparison is between the one-field FDM and implicit IFEM using the same time step $\Delta t=1.0\times 10^{-4}$.}  
	\label{leaflet_defomration_parameter2}
\end{figure}

\begin{figure}[h!]
	\begin{minipage}[t]{1\linewidth}
		\centering
		\includegraphics[width=5in,angle=0]{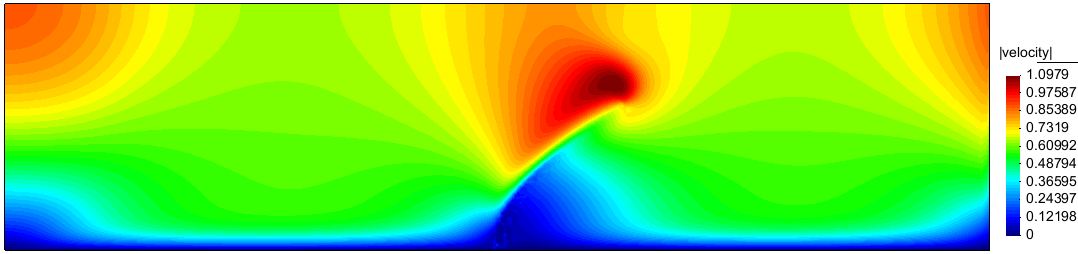}
		\captionsetup{justification=centering}
		\caption* {\scriptsize (a) $t/T=0.1$} 
	\end{minipage}
	\begin{minipage}[t]{1\linewidth}
		\centering
		\includegraphics[width=5in,angle=0]{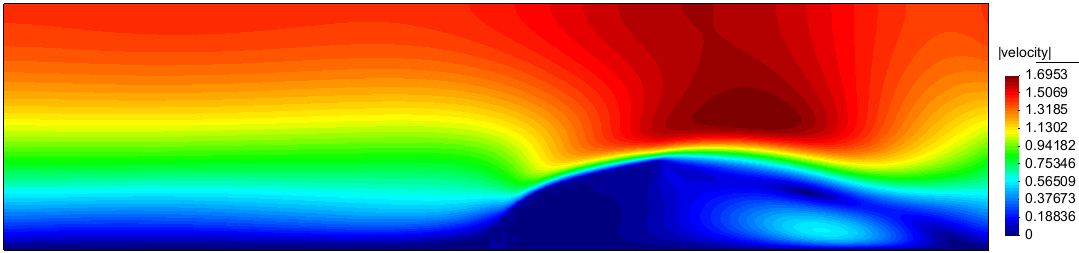}
		\captionsetup{justification=centering}
		\caption* {\scriptsize (b) $t/T=0.3$.} 
	\end{minipage} 
	\begin{minipage}[t]{1\linewidth}
		\centering		
		\includegraphics[width=5in,angle=0]{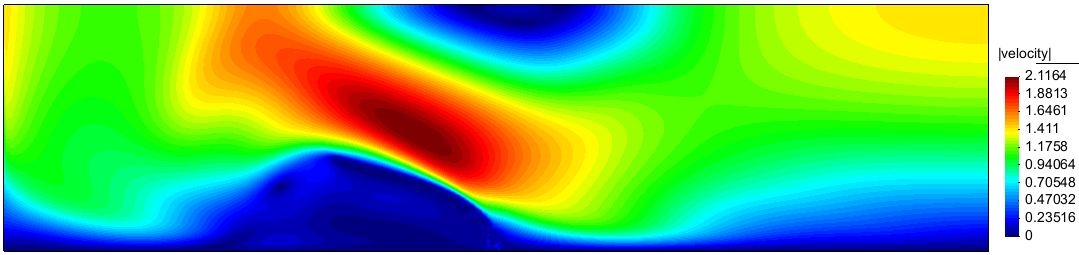}
		\captionsetup{justification=centering}
		\caption* {\scriptsize (c) $t/T=0.7$} 
	\end{minipage}
	\begin{minipage}[t]{1\linewidth}
		\centering
		\includegraphics[width=5in,angle=0]{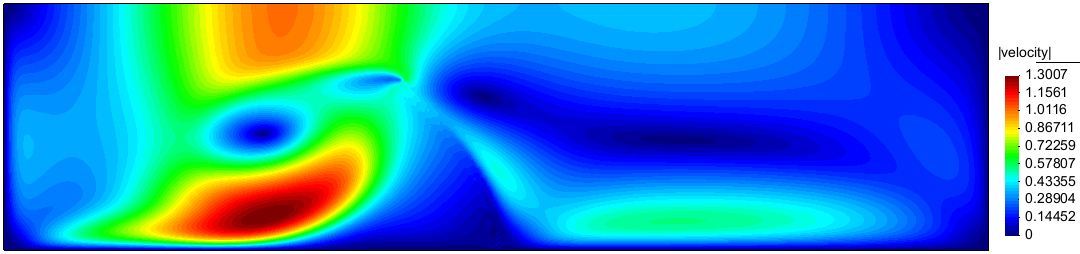}
		\captionsetup{justification=centering}
		\caption* {\scriptsize (d) $t/T=1.0$.} 
	\end{minipage} 	
	\captionsetup{justification=centering}
	\caption {\scriptsize Velocity norm in the background domain for Parameter set 3 using the one-field FDM.}  
	\label{leaflet_defomration_parameter3_fliud}
\end{figure}

\begin{figure}[h!]
	\begin{minipage}[t]{0.5\linewidth}
		\centering
		\includegraphics[width=2.2in,angle=0]{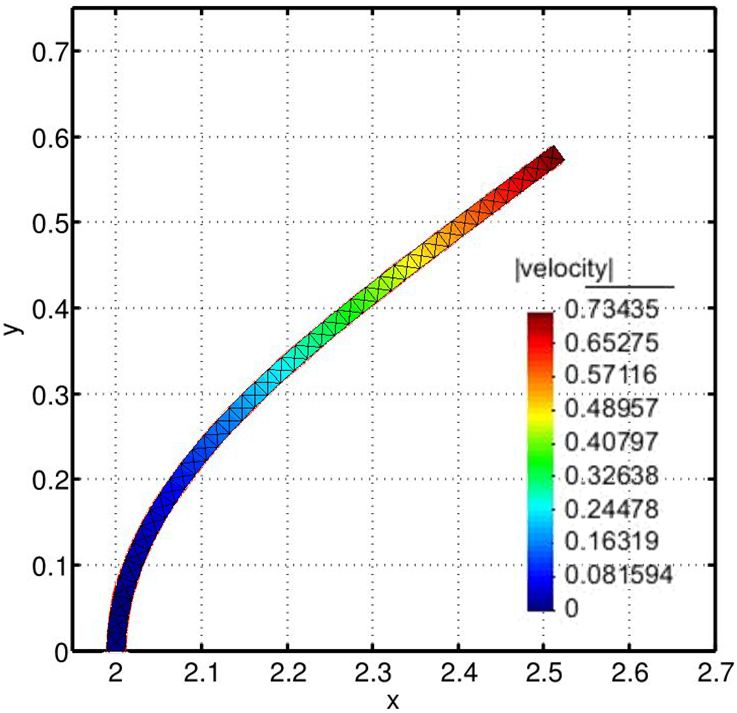}
		\captionsetup{justification=centering}
		\caption* {\scriptsize (a) $t/T=0.1$} 
	\end{minipage}
	\begin{minipage}[t]{0.5\linewidth}
		\centering
		\includegraphics[width=2.2in,angle=0]{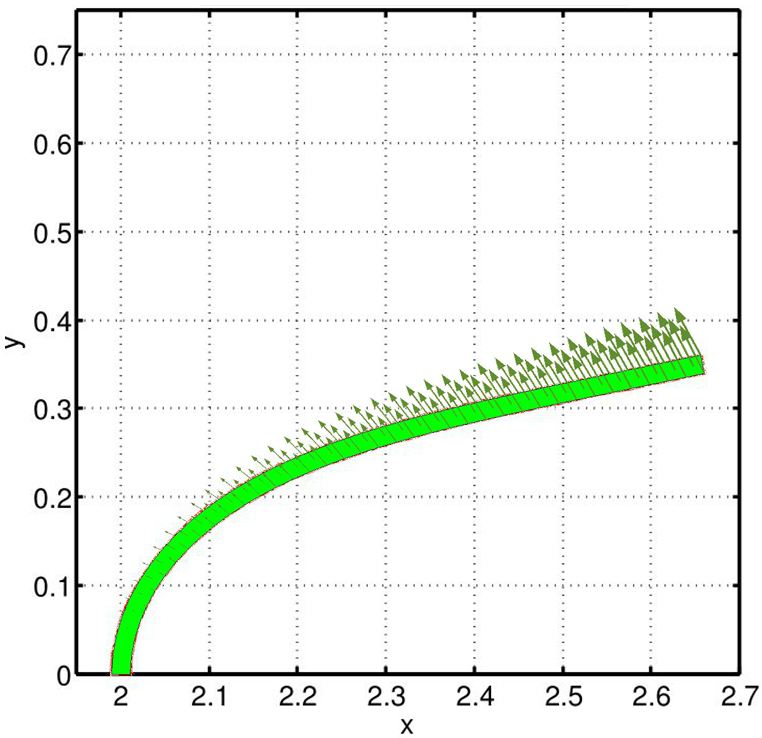}
		\captionsetup{justification=centering}
		\caption* {\scriptsize (b) $t/T=0.3$.} 
	\end{minipage} 
	\begin{minipage}[t]{0.5\linewidth}
		\centering
		\includegraphics[width=2.2in,angle=0]{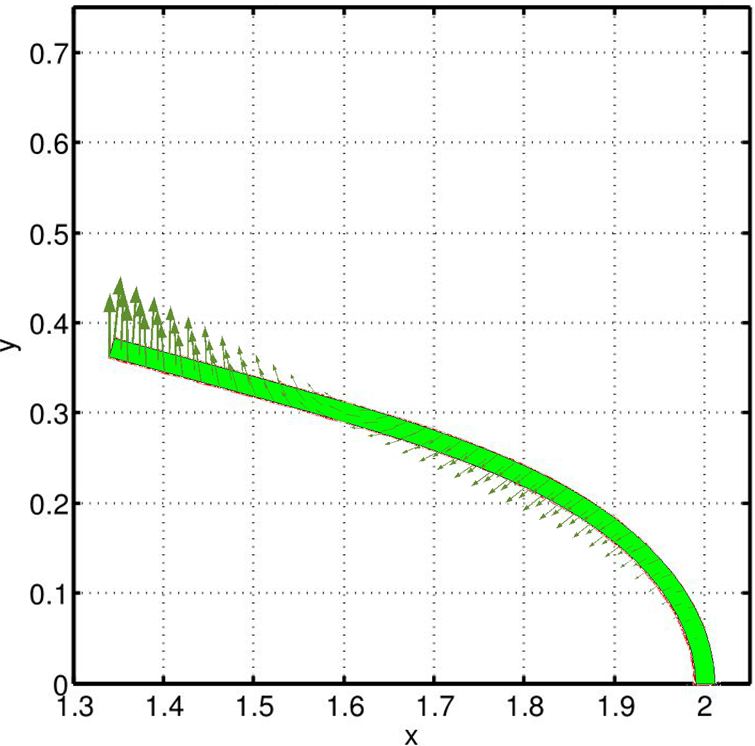}
		\captionsetup{justification=centering}
		\caption* {\scriptsize (c) $t/T=0.7$} 
	\end{minipage}
	\begin{minipage}[t]{0.5\linewidth}
		\centering
		\includegraphics[width=2.2in,angle=0]{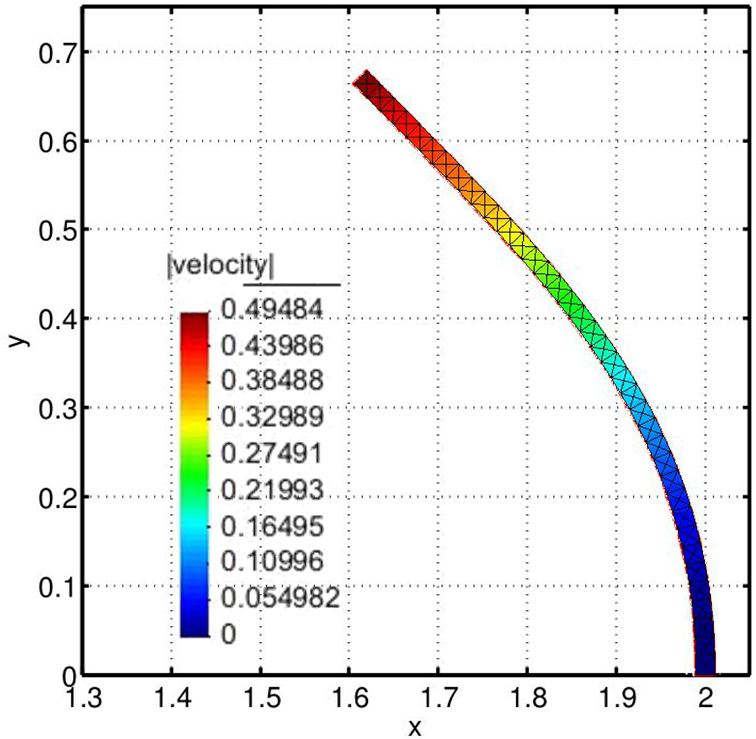}
		\captionsetup{justification=centering}
		\caption* {\scriptsize (d) $t/T=1.0$.} 
	\end{minipage} 	
	\captionsetup{justification=centering}
	\caption {\scriptsize Leaflet deformation at different times for Parameter set 3 using the one-field FDM.}  
	\label{leaflet_defomration_parameter3_solid}
\end{figure}

We then test a case with different density between fluid and solid: $\rho^r=1.2$. We use the same time step $\Delta t=1.0\times 10^{-4}$ for both the one-field FDM and the implicit IFEM, and their results, observed from Figure \ref{leaflet_defomration_parameter2}, are very close. The case of $\rho^r=2$ has also been tested, but our nonlinear implicit IFEM solver cannot converge at the first time step for any time step size. For completeness, we present the results using the proposed one-field FDM in Figures \ref{leaflet_defomration_parameter3_fliud} and \ref{leaflet_defomration_parameter3_solid}.

%%%%%%%%%%%%%%%%%%%%%%%%%%%%%%%%%%%%%%%%%%
\section{Conclusions}
\label{sec_conclusions}
In this article we have illustrated the relationship between the recently proposed one-field FDM and both the explicit (explicitly expressing the force term) and implicit (implicitly expressing the force term) IFEM. This is facilitated through the use of a particular operator splitting scheme. Furthermore, we show that the one-field FDM produces the same accuracy of results as the implicit IFEM, but requires no iteration, whilst it significantly improves upon the classical IFEM at very little additional computational complexity. The one-field FDM is shown to naturally deal with the case of different densities and/or viscosities of the fluid and solid. Therefore, whilst we may view the scheme of \cite{Wang_2017} as a fictitious domain method, it is also legitimate to consider it to be an alternative, highly efficient and robust, approximate solution strategy for the fully implicit IFEM methods of \cite{wang2006immersed, wang2007iterative, wang2009computational}. Even in the simple case of the same density and viscosity, where the explicit IFEM is known to be successful, we find that the additional terms added in the formulation of the one-field FDM have a helpful stabilizing effect, such that a larger time step can be adopted compared with the explicit IFEM.

\appendix
\section{Extension to the compressible neo-Hookean solid model}
In this section, we extend the incompressible neo-Hookean solid model to a compressible case, in which case the constitutive equation can be expressed as \cite{Boffiarchive}:
\begin{equation} \label{constitutive_solid_com}
{\bm\sigma}={\bm\sigma}^s=c_1J^{-1}\left({\bf F}{\bf F}^T-{\bf I}\right)+\mu^s{\rm D}{\bf u}^s-J^{-1/(1-2\nu)}{\bf I},
\end{equation}
where $\nu$ is the Poisson's ratio. For a compressible solid, the continuity equation can simply be expressed as:
\begin{equation}
J\rho^s=\rho_0^s,
\end{equation}
where $\rho_0^s$ is the initial solid density. Then the corresponding weak forms (\ref{weak_form1}) and (\ref{weak_form2}) can be expressed as \cite{Boffiarchive}:
\begin{equation}\label{weak_form1_com_neo}
\begin{split}
&\rho^f\int_{\Omega}\frac{\partial{\bf u}}{\partial t} \cdot{\bf v}
+\rho^f \int_{\Omega}\left({\bf u}\cdot\nabla\right){\bf u}\cdot{\bf v}
+\frac{\mu^f}{2}\int_{\Omega}{\rm D}{\bf u}:{\rm D}{\bf v}
-\int_{\Omega}p\nabla \cdot {\bf v} \\
+&\rho^{\delta}\int_{\Omega_t^s}\frac{\mathfrak{d}{\bf u}}{\mathfrak{d} t}\cdot{\bf v}
+\frac{\mu^\delta}{2}\int_{\Omega_t^s}{\rm D}{\bf u}:{\rm D}{\bf v}
+c_1\int_{\Omega_t^s}\left(J^{-1}{{\bf F}{\bf F}^T-\boxed{J^{-1/(1-2\nu)}{\bf I}}}\right):\nabla{\bf v}  \\
&+\boxed{\int_{\Omega_t^s}p\nabla \cdot {\bf v} }
=\int_{\Omega}\rho^f{\bf g}\cdot{\bf v}
+\int_{\Omega_t^s}\rho^\delta{\bf g}\cdot{\bf v}
+\int_{\Gamma_N}\bar{\bf h}\cdot{\bf v},
\end{split}
\end{equation}
with $\rho^\delta=\rho_0^s/J-\rho^f$, and
\begin{equation}\label{weak_form2_com_neo}
-\int_{\Omega} q\nabla \cdot {\bf u}+\boxed{\int_{\Omega_t^s} q\nabla \cdot {\bf u}}+\boxed{\frac{1}{\kappa}\int_{\Omega_t^s} qp}=0.
\end{equation}	

The boxed terms in the above equations indicate the differences compared with equation (\ref{weak_form1}) and (\ref{weak_form2}). For a compressible solid model, the incompressibility equation (\ref{continuity_equation}) ($\nabla\cdot{\bf u}=0$) only holds in the fluid domain $\Omega_t^f$. Therefore we cannot solve it in the whole domain $\Omega$ using a fictitious domain method, because this never matches the velocity of a compressible solid ($\nabla\cdot{\bf u}\ne 0$). The pressure computed in the solid domain $\left.p\right|_{\Omega_t^s}$ is meaningless, which is weakly imposed to be zero in (\ref{weak_form2_com_neo}) with $\kappa$ playing the role of a bulk modulus \cite{Boffiarchive}.

In order to implement the one-field FDM approach, after time discretization one could update the solid stress as described in Section \ref{one_field_fdm}. Alternatively, one could also update the deformation tensor as follows.
\begin{equation}\label{update_f}
\begin{split}
&\int_{\Omega_t^s}J^{-1}{\bf F}_{n+1}{\bf F}_{n+1}^T:\nabla{\bf v}
=\int_{\Omega_{\bf X}^s}{\bf F}_{n+1}:\nabla_{\bf X}{\bf v} \\
&=\int_{\Omega_{\bf X}^s}\left({\bf F}_n+\Delta t\nabla_{\bf X}{\bf u}_{n+1}\right):\nabla_{\bf X}{\bf v}.
\end{split}
\end{equation}
Using (\ref{update_f}), equation (\ref{weak_form1_com_neo}) and (\ref{weak_form2_com_neo}) may be solved implicitly, which can also use the operator spitting scheme introduced in Section \ref{splitting_scheme}. See \cite{yongxing2018} for more details about the implicit solver and different types of explicit splitting schemes.

\section{Extension to the compressible Saint Venant-Kirchhoff solid model}
The constitutive equation of the Saint Venant-Kirchhoff solid model can be expressed as \cite{Bazilevs_2010}:
\begin{equation}\label{expression_of_s}
{\bf S}\left({\bf E}\right)=2\mu{\bf E}+\lambda tr\left({\bf E}\right){\bf I},
\end{equation}
where 
\begin{equation}\label{green_strain}
{\bf E}=\frac{1}{2}\left({\bf F}^T{\bf F}-{\bf I}\right)
\end{equation}
is the Lagrangian Green strain, $\mu$ and $\lambda$ are the Lam$\acute{\rm e}$ constants. Then the corresponding weak forms (\ref{weak_form1}) and (\ref{weak_form2}) can be expressed as:
\begin{equation}\label{weak_form1_com}
\begin{split}
&\rho^f\int_{\Omega}\frac{\partial{\bf u}}{\partial t} \cdot{\bf v}
+\rho^f \int_{\Omega}\left({\bf u}\cdot\nabla\right){\bf u}\cdot{\bf v}
+\frac{\mu^f}{2}\int_{\Omega}{\rm D}{\bf u}:{\rm D}{\bf v}
-\int_{\Omega}p\nabla \cdot {\bf v} \\
+&\rho^{\delta}\int_{\Omega_t^s}\frac{\mathfrak{d}{\bf u}}{\mathfrak{d} t}\cdot{\bf v}
+\frac{\mu^\delta}{2}\int_{\Omega_t^s}{\rm D}{\bf u}:{\rm D}{\bf v}
+\boxed{\frac{1}{2}\int_{\Omega_{\bf X}^s}{\bf S}:\delta{\bf E}} \\
&+\boxed{\int_{\Omega_t^s}p\nabla \cdot {\bf v} }
=\int_{\Omega}\rho^f{\bf g}\cdot{\bf v}
+\int_{\Omega_t^s}\rho^\delta{\bf g}\cdot{\bf v}
+\int_{\Gamma_N}\bar{\bf h}\cdot{\bf v}
\end{split}
\end{equation}
with $\delta{\bf E}={\bf F}^T\left(\nabla_{\bf X}{\bf v}\right)+\left(\nabla_{\bf X}^T{\bf v}\right){\bf F}$ and $\rho^\delta=\rho_0^s/J-\rho^f$, and
\begin{equation}\label{weak_form2_com}
-\int_{\Omega} q\nabla \cdot {\bf u}+\boxed{\int_{\Omega_t^s} q\nabla \cdot {\bf u}}
+\boxed{\frac{1}{\kappa}\int_{\Omega_t^s} qp}=0.
\end{equation}	

The term ${\bf S}:\delta{\bf E}$ may be linearized at a given displacement $\tilde{\bf d}$ as follows:
\begin{equation}\label{f_linearization}
{\bf S}:\delta{\bf E}
\approx{\bf S}\left(\tilde{\bf E}\right)
:\delta\tilde{\bf E}
-\frac{1}{2}{\bf S}\left(\nabla_{\bf X}^T\tilde{\bf d}\nabla_{\bf X}\tilde{\bf d}\right):
\delta\tilde{\bf E},
\end{equation}
where
\begin{equation}
\tilde{\bf E}=\frac{1}{2}\left({\bf\rm D}_{\bf X}{\bf d}+\nabla_{\bf X}^T\tilde{\bf d}\nabla_{\bf X}{\bf d}+\nabla_{\bf X}^T{\bf d}\nabla_{\bf X}\tilde{\bf d}\right)
\end{equation}
and
\begin{equation}
\delta\tilde{\bf E}
=\frac{1}{2}\left({\bf\rm D}_{\bf X}{\bf v}+\nabla_{\bf X}^T\tilde{\bf d}\nabla_{\bf X}{\bf v}+\nabla_{\bf X}^T{\bf v}\nabla_{\bf X}\tilde{\bf d}\right).
\end{equation} 
As with the previous implementation of the one-field FDM, one may update the solid displacement ${\bf d}$ after time discretizaiton as follows:
\begin{equation}
{\bf d}_{n+1}={\bf d}_n+\Delta t{\bf u}_{n+1}.
\end{equation}
One still can use the operator splitting scheme by choosing $\tilde{\bf d}={\bf d}_n$, although we omit the full details here.

\end{document}